\documentclass[twoside,11pt]{article}
\title{} \author{} \date{}
\usepackage{latexsym,amssymb,times}%,showkeys}
\input amssym.def
\newtheorem{te}{Theorem}[section]

\newtheorem{fac}[te]{Fact}
\newtheorem{lem}[te]{Lemma}
\newtheorem{cla}[te]{Claim}

\newtheorem{ex}[te]{Example}

\def\dok{\noindent{\bf Proof. }}
\def\kdok{\hfill $\Box$ \par \vspace*{2mm} }
\def\f{\varphi}
\def\p{\psi}
\def\o{\omega}
\def\k{\kappa}
\def\l{\lambda}
\def\r{\rho}
\def\s{\sigma}
\def\t{\tau}
\def\P{{\mathbb P}}
\def\Q{{\mathbb Q}}
\def\B{{\mathbb B}}
\def\N{{\mathbb N}}
\def\X{{\mathbb X}}
\def\Y{{\mathbb Y}}
\def\C{{\mathbb C}}
\def\R{{\mathbb R}}
\def\CD{{\mathcal D}}
\def\CC{{\mathcal C}}
\def\M{{\mathcal M}}
\def\down{\!\downarrow}
\def\up{\!\uparrow}
\def\la{\langle}
\def\ra{\rangle}
\def\id{\mathop{\rm id}\nolimits}
\def\Borel{\mathop{\rm Borel}\nolimits}
\def\Scatt{\mathop{\rm Scatt}\nolimits}
\def\Emb{\mathop{\rm Emb}\nolimits}
\def\Iso{\mathop{\rm Iso}\nolimits}
\def\Aut{\mathop{\rm Aut}\nolimits}
\def\Fin{\mathop{\rm Fin}\nolimits}
\def\sm{\mathop{\rm sm}\nolimits}
\def\sq{\mathop{\rm sq}\nolimits}
\def\ar{\mathop{\rm ar}\nolimits}
\def\Int{\mathop{\rm Int}\nolimits}
\def\Mod{\mathop{\rm Mod}\nolimits}
\def\Col{\mathop{\rm Col}\nolimits}
\def\ro{\mathop{\rm ro}\nolimits}
\def\ge{\mathop{\rm ge}\nolimits}
\begin{document}
\thispagestyle{plain}
\begin{center}
           {\large \bf {\uppercase{Different Similarities}}}
\end{center}
\begin{center}
{\bf Milo\v s S.\ Kurili\'c\footnote{Department of Mathematics and Informatics, University of Novi Sad,
         Trg Dositeja Obradovi\'ca 4, 21000 Novi Sad, Serbia.
                                     e-mail: milos@dmi.uns.ac.rs}}
\end{center}

\begin{abstract}
\noindent
We establish the hierarchy among twelve equivalence relations (similarities) on the class of relational structures:
the equality, the isomorphism, the equimorphism, the full relation, four similarities of structures induced by similarities of their
self-embedding monoids and intersections of these equivalence relations. In particular, fixing a language $L$ and a cardinal $\k$,
we consider the interplay between the restrictions of these similarities to the class $\Mod _L (\k )$ of all $L$-structures of size $\k$.
It turns out that, concerning the number of different similarities and the shape of the corresponding Hasse diagram,
the class of all structures naturally splits into three parts: finite structures,
infinite structures of unary languages, and  infinite structures of non-unary languages (where all these similarities are different).

%\vspace{1mm}\\

\end{abstract}
\section{Introduction}\label{S1}
If $\X$ is a relational structure, $\Emb (\X )$ the monoid of its self-embeddings and  $\P (\X )= \{ f[X] : f \in \Emb (\X )\} $
the set of copies of $\X$ inside $\X$, then the poset $\la \P (\X ),\subset \ra$ (isomorphic to the inverse of the right
Green's order on $\Emb (\X )$) contains a certain information about $\X$ and the equality $\P (\X )=\P (\Y )$ defines an
equivalence relation on the class of all relational structures. Writing $\P (\X )$ instead of $\la \P (\X ), \subset \ra$,
some coarser classifications of structures are obtained if the equality is replaced
by the following weaker conditions: $\P (\X )\cong \P (\Y )$ (implied  by $\Emb (\X )\cong \Emb  (\Y )$),
$\sq \P (\X )\cong \sq \P (\Y )$ (where $\sq \P$ denotes the separative quotient of a poset $\P $), and $\P (\X )\equiv \P (\Y )$
(the forcing equivalence of posets of copies). Concerning the last (and the coarsest non-trivial) similarity relation we note that
the forcing related properties of posets of copies was investigated for countable structures
in general in \cite{Ktow},
for equivalence relations and similar structures in \cite{Kemb},
for ordinals in \cite{Kurord},
for scattered and non-scattered linear orders in \cite{Kur1} and \cite{KurTod},
and for several ultrahomogeneous structures in \cite{Kstr},\cite{KurTod},\cite{KurTod1}, and \cite{KurTod2}.

In this paper we investigate the interplay between the four similarity relations mentioned above and the similarities defined by
the conditions $\X = \Y$, $\X \cong \Y$, and $\X \rightleftarrows \Y$ (equimorphism, bi-embeddability).

\begin{figure}[htb]%\label{F4002}
\unitlength 0.9mm%1mm %0.7mm % = .854pt
\linethickness{0.5pt}
\ifx\plotpoint\undefined\newsavebox{\plotpoint}\fi % GNUPLOT compatibility
\begin{picture}(140,120)(0,0)

%----------------------------- linije --------------------------

\put(55,10){\line(0,1){15}}%1
\put(55,25){\line(0,1){15}}%2
\put(55,25){\line(2,1){30}}%3
\put(55,40){\line(-2,1){30}}%4
\put(55,40){\line(2,1){30}}%5
\put(85,40){\line(0,1){15}}%6
\put(25,55){\line(2,1){30}}%7
\put(85,55){\line(-2,1){30}}%8
\put(85,55){\line(0,1){15}}%9
\put(55,70){\line(0,1){15}}%10
\put(85,70){\line(-2,1){30}}%11
\put(85,70){\line(2,1){30}}%12
\put(55,85){\line(2,1){30}}%13
\put(115,85){\line(-2,1){30}}%14
\put(85,100){\line(0,1){15}}%15

%----------------------------- tacke --------------------------

\put(55,10){\circle*{1}}%1
\put(55,25){\circle*{1}}%
\put(55,40){\circle*{1}}%
\put(85,40){\circle*{1}}%
\put(25,55){\circle*{1}}%
\put(85,55){\circle*{1}}%
\put(55,70){\circle*{1}}%
\put(85,70){\circle*{1}}%
\put(55,85){\circle*{1}}%
\put(115,85){\circle*{1}}%
\put(85,100){\circle*{1}}%
\put(85,115){\circle*{1}}%

%----------------------------- tekst -----------------------
%\small
\footnotesize
%\scriptsize
%\tiny

\put(57,5){\makebox(0,0)[cc]{$\X =\Y \land L_1=L_2$}}%A
\put(13,55){\makebox(0,0)[cc]{$\P (\X )= \P (\Y )$}}%B
\put(93,40){\makebox(0,0)[cc]{$\X \cong \Y $}}%C
\put(42,70){\makebox(0,0)[cc]{$\P (\X )\cong \P (\Y )$}}%D
\put(123,85){\makebox(0,0)[cc]{$\X \leftrightarrows \Y$}}%E
\put(39,85){\makebox(0,0)[cc]{$\sq \P (\X )\cong \sq \P (\Y )$}}%F
\put(65,100){\makebox(0,0)[cc]{$\ro \sq  \P (\X )\cong \ro  \sq \P (\Y )$}}
\put(100,100){\makebox(0,0)[cc]{$\Leftrightarrow \P (\X )\equiv \P (\Y )$}}%G
\put(75,23){\makebox(0,0)[cc]{$\P (\X )= \P (\Y )\land \X \cong \Y $}}%H
\put(35,38){\makebox(0,0)[cc]{$\P (\X )= \P (\Y ) \land \X \leftrightarrows \Y$}}%H
\put(105,55){\makebox(0,0)[cc]{$\P (\X )\cong \P (\Y )\land \X \leftrightarrows \Y$}}%H
\put(109,70){\makebox(0,0)[cc]{$\sq \P (\X )\cong \sq \P (\Y )\land \X \leftrightarrows \Y$}}%H
\put(99,115){\makebox(0,0)[cc]{the full relation}}%A

\tiny
\put(58,17){\makebox(0,0)[cc]{$a$}}%
\put(58,32){\makebox(0,0)[cc]{$b$}}%
\put(75,32){\makebox(0,0)[cc]{$c$}}%
\put(47,47){\makebox(0,0)[cc]{$d$}}%
\put(63,47){\makebox(0,0)[cc]{$e$}}%
\put(88,47){\makebox(0,0)[cc]{$f$}}%
\put(47,63){\makebox(0,0)[cc]{$g$}}%
\put(63,63){\makebox(0,0)[cc]{$h$}}%
\put(88,62){\makebox(0,0)[cc]{$i$}}%
\put(51,77){\makebox(0,0)[cc]{$j$}}%
\put(77,77){\makebox(0,0)[cc]{$k$}}%
\put(93,77){\makebox(0,0)[cc]{$l$}}%
\put(77,93){\makebox(0,0)[cc]{$m$}}%
\put(93,93){\makebox(0,0)[cc]{$n$}}%
\put(88,107){\makebox(0,0)[cc]{$o$}}%
\end{picture}
\caption{The hierarchy of similarities between relational structures}\label{F4002}
\end{figure}
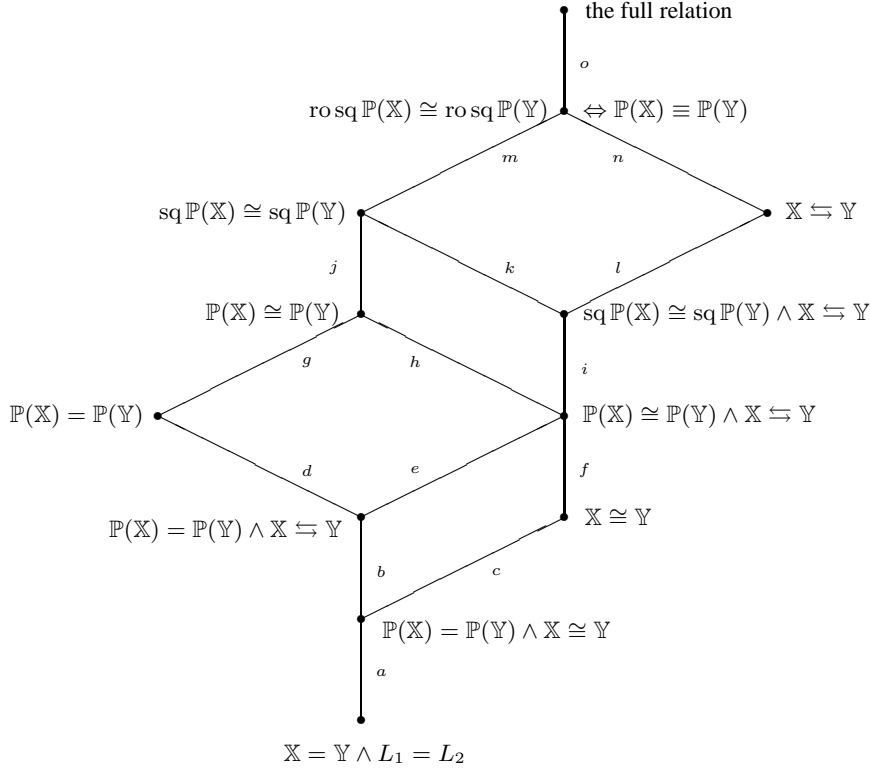

In Section \ref{S2} we establish the hierarchy displayed in Figure \ref{F4002}, which, more precisely,
contains the implications between the similarities on the class of {\it pairs} $\la \X, L \ra$, where $L$ is a language and $\X$ an $L$-structure.
(The language must be included in the game because, otherwise, since the structure $\X =\la \o , \la \emptyset\ra\ra$ can be regarded as an $L$-structure
for each language $L$ of size 1, it is not clear what $\X \cong \Y$ means). So, the conditions displayed in the diagram define when
the pairs $\la \X ,L_1\ra$ and $\la \Y ,L_2\ra$  are similar (clearly, the equality $L_1=L_2$ follows from $\X \cong \Y$ and $\X \rightleftarrows\Y$
and we omit it). Thus, for example, line {\it n} denotes the statement that equimorphic structures have
forcing-equivalent posets of copies.

In Section \ref{S3} we fix a language $L$ and a set $X$ and restrict our analysis to the class $\Mod _L (X)$ of $L$-structures with
the domain $X$. It turns out that for a non-unary language $L$ and infinite set $X$ in the diagram from Figure \ref{F4002}
 restricted to the class $\Mod _L (\k )$ all the implications $a$ - $o$ are proper and there are no new implications
(except the ones following from transitivity). On the other hand, for finite structures or unary languages the diagram
collapses significantly.

A few words on notation.
Let $L=\la R_i :i\in I \ra$ be a relational language, where $\ar _L (R_i)=n_i\in \N$, $i\in I$ and let $X$ be a non empty set.
If $\X = \la X, \la \r _i :i\in I \ra \ra$ is an $L$-structure
and  $\emptyset\neq A\subset X$, then $\la A, \la \r _i \upharpoonright A :i\in I \ra \ra$ is a
{\it substructure} of $\X$, where $\r_i \upharpoonright A = \r _i \cap A^{n_i}$, for $i\in I$.
If $\Y =\la Y, \la \s _i : i\in I\ra  \ra$ is an $L$-structure too,
a mapping $f:X \rightarrow Y$ is an
{\it embedding} (we write $f:\X \hookrightarrow  \Y$) iff $f$ is an injection and
for all $i\in I$ and $x_1, \dots   x_{n _i} \in X$ we have $\la x_1, \dots   ,x_{n _i}\ra \in \r _i
\Leftrightarrow \la f(x_1), \dots  ,f(x_{n _i})\ra \in \s _i$.
%Then we write $\X \hookrightarrow \Y$.
Let $\Emb (\X , \Y )$ denote the set of such embeddings and
$\P (\X , \Y  )  =  \{ B\subset Y : \la B, \la \s _i \upharpoonright B :i\in I \ra \ra \cong \X \}
                =  \{ f[X] : f \in \Emb (\X , \Y )\} .$
In particular, $\Emb (\X ) := \Emb (\X, \X)$ and
$\P (\X )  :=  \{ f[X] : f \in \Emb (\X )\} =  \{ A\subset X : \la A, \la \r _i \upharpoonright A :i\in I \ra \ra \cong \X\}$.
If  $f\in \Emb (\X , \Y )$ is a surjection, it is an
{\it isomorphism}, we write $f\in \Iso (\X , \Y ) $, and the structures $\X$ and $\Y$ are
{\it isomorphic}, in notation $\X \cong \Y$.
If, in particular, $ \Y = \X $, then $f$ is called an
{\it automorphism} of the structure $\X$ and $\Aut (\X )$ denotes the set of all such mappings.
Structures $\X$ and $\Y$ are called
{\it equimorphic}, in notation $\X\rightleftarrows\Y$, iff $\X \hookrightarrow \Y$ and $\Y \hookrightarrow \X$.

The {\it right Green's pre-order} $\preceq ^R$ on the monoid $\la \Emb \X , \circ , \id _X \ra$
is defined by: $f\preceq ^R g$ iff $f\circ h=g$, for some $h\in \Emb \X $.
The {\it right Green's equivalence relation}  $\approx ^R$  on $\Emb \X$, given by: $f\approx^R g$ iff $f\preceq ^R g$ and $g\preceq ^R f$,
determines the antisymmetric quotient $\la \Emb \X /\!\approx ^R, \preceq ^R\ra$, the {\it right Green's order}.
It is easy to check that  $\la \Emb \X /\!\!\approx ^R, \preceq ^R\ra \cong \la \P (\X ), \supset \ra$ so the results of this paper
can be regarded as statements about transformation semigroups.

A partial order $\P =\la P , \leq \ra$ is called
{\it separative} iff for each $p,q\in P$ satisfying $p\not\leq q$ there is $r\in P$ such that $r\leq p$ and $r\perp q$.
The {\it separative modification} of $\P $ is the pre-order $\sm \P =\la P , \leq ^*\ra$, where
$p\leq ^* q$ iff $\forall r\leq p \; \exists s \leq r \; s\leq q $.
The {\it separative quotient} of $\P $
is the separative partial order $\sq \P = \langle P /\!\! =^* , \trianglelefteq \rangle$, where
$p = ^* q \Leftrightarrow p \leq ^* q \land q \leq ^* p\;$  and $\;[p] \trianglelefteq [q] \Leftrightarrow p \leq ^* q $.
If $\P$ is a separative partial order, by $\ro \P$ we will denote the {\it Boolean completion} of $\P$.
For a pre-order $\P$ let $\ge (\P)= \{ V_\P [G]: G$ is a $\P$-generic filter over $V \}$. Two pre-orders $\P $ and $\Q$ are said to be
{\it forcing equivalent}, in notation $\P \equiv \Q $, iff  $\ge (\P)=\ge (\Q)$.
\begin{fac}   \label{T4042}
Let $\P , \Q $ and $\P _i$, $i\in I$, be partial orderings. Then

(a) $\P \cong \Q \Rightarrow \sm \P \cong  \sm \Q \Rightarrow \sq \P \cong  \sq \Q \Rightarrow \ro \sq \P \cong  \ro \sq \Q \Rightarrow \P\equiv \Q$;

(b) $\P \equiv \sm \P \equiv \sq \P \equiv  (\ro \sq \P)^+ $;

(c) $\sq (\prod _{i\in I}\P _i) \cong \prod _{i\in I}\sq \P _i$.
\end{fac}
\section{Implications}\label{S2}
In this section we establish the implications $a$ - $o$ from Figure \ref{F4002}. In Section \ref{S3} we will show that, regarding
the class of all relational structures, there are no new implications in Figure \ref{F4002} (except the ones which follow from the transitivity).
First, the implications $a$, $b$, $c$, $d$, $e$, $g$, $h$, $k$, $l$, and $o$ are evident, while $i$, $j$ and $m$ follow from Fact \ref{T4042}(a).
In the sequel we prove the equivalence $\ro \sq  \P (\X )\cong \ro  \sq \P (\Y )\Leftrightarrow \P (\X )\equiv \P (\Y )$ and the implications $f$ and $n$
(see Theorems \ref{T6032}  and \ref{T6004}).
\subsection{Intermezzo: the homogeneity of Boolean completions}
Here we prove that the Boolean completion of the poset of copies of a relational structure is a
homogeneous Boolean algebra.
We recall that a partial order $\P =\la P ,\leq \ra$ is called
{\it homogeneous} iff  it has a largest element and $\P \cong p\down $, for each $p\in P$ and that
a Boolean algebra $\B $ is called a
{\it homogeneous Boolean algebra} iff  $\B \cong b\down $, for each $b \in \B ^+$.
It is known  that the Boolean completion of a separative homogeneous partial order $\P$ is a homogeneous Boolean algebra (see \cite{Kopp1}, p.\ 181)
and, by Theorem 2.2 of \cite{Ktow}, the posets of the form $\P (\X )$ are homogeneous but it is easy to see that they are not separative in most of the cases.
So, in order to prove that the Boolean completions $\ro \sq \P (\X )$ are homogeneous algebras, we show that in the theorem mentioned above
the  separativity of $\P$ can be omitted
and that the assumption of homogeneity can be relaxed. Namely, defining a partial order $\P $ to be {\it quasi homogeneous} iff  for each $p\in P$
there is a dense subset of $\P$ isomorphic to a dense subset of $p\down $,  we have the following generalization.
\begin{te}   \label{T6030}
The Boolean completion of a quasi homogeneous partial order $\P$ is a homogeneous Boolean algebra.
\end{te}
\dok The statement is a consequence of the following two claims. Namely, if $\P$ is a quasi homogeneous partial ordering, then, by Claim \ref{T6001},
$\sq \P$ is a separative quasi homogeneous partial order and,  by Claim \ref{T6031}, the algebra $\ro \sq \P$ is homogeneous.
\begin{cla}   \label{T6001}
The separative quotient of a quasi homogeneous partial order is quasi homogeneous.
\end{cla}
\dok
Let $\P =\la P\leq\ra$ be a quasi homogeneous partial order, $\sq \P  =\langle P /\!\! =^* , \trianglelefteq \rangle $ and $p\in P$.
Let $D$ be a dense subset of $\P$ and $f: \la D , \leq \ra  \rightarrow \la p\down ,\leq \ra$ an embedding such that $f[D]$ is a dense subset of $p\down$.
First we prove that
\begin{equation}\label{EQ6048}
\forall q,r \in D \;\; (q\leq ^* r \Leftrightarrow f(q)\leq ^* f(r)).
\end{equation}
Let $q,r \in D$. If $q\leq ^* r$, then  each $s\leq q$ is compatible with $r$ and we prove that  each $u\leq f(q)$ is compatible with $f(r)$.
If $u\leq f(q)$, then $u\leq p$ and, since $f[D]$ is dense in $p\down$,
there is $s\in D$ such that $f(s)\leq u$.
Since $f$ is an embedding and $f(s)\leq f(q)$ we have $s\leq q$
and, since $q\leq ^* r$, there is $t\leq s,r$, and, moreover there is $t'\in D$ such that $t'\leq t$
which implies $f(t')\leq f(s)\leq u$ and $f(t')\leq f(r)$. Thus $u\not\perp f(r)$.

Assuming that $f(q)\leq ^* f(r)$ and $s\leq q$ we show that $s\not\perp r$.
If $s\leq q$ and $s'\in D$, where $s'\leq s$,
then $f(s')\leq f(q)$
and, since $f(q)\leq ^* f(r)$, there is $v\leq f(s'),f(r)\leq p$.
Since $f[D]$ is dense in $p\down$, there is $t\in D$ such that $f(t)\leq v$.
Since $f$ is an embedding we have $t\leq s',r$ and, hence, $s\not\perp r$. So (\ref{EQ6048}) is true.

It is evident that the set $\CD :=\{ [q]:q\in D \}$ is a dense suborder of the partial order $\langle P /\!\! =^* , \trianglelefteq \rangle$
and we prove that the mapping
$$
F: \langle \CD , \trianglelefteq \rangle\rightarrow \langle [p]\down  , \trianglelefteq \rangle ,
$$
given by $F([q])=[f(q)]$, is an embedding.
First, for $q,r\in D$ by (\ref{EQ6048}) we have $[q]=[r]$ iff $q=^* r$ iff $q\leq ^* r \land r\leq ^* q$ iff
$f(q)\leq ^* f(r)\land f(r)\leq ^* f(q)$ iff $f(q)= ^* f(r)$
iff $[f(q)]= [f(r)] $ iff $F([q])=F([r])$ and, thus, $F$ is a well defined injection.
Second, for $q\in D$ we have $f(q)\leq p$, which implies $f(q)\leq ^* p$ and, hence, $[f(q)]\trianglelefteq [p]$, that is $F([q])\in [p]\down$.
Thus $F[\CD]\subset [p]\down$.
Finally, by (\ref{EQ6048}), for $q,r\in D$ we have $[q]\trianglelefteq [r]$ iff $q\leq ^* r$ iff $f(q)\leq ^* f(r)$ iff $[f(q)]\trianglelefteq [f(r)]$
iff $F([q]) \trianglelefteq F([r])$ and, thus, $F$ is a strong homomorphism.

Now we prove that $F[ \CD]$ is a dense set in the poset $\langle [p]\down  , \trianglelefteq \rangle$.
If $[q]\trianglelefteq [p]$, then there is $s\leq p,q$ and,
since $f[D]$ is dense in $p\down$, there is $u\in D$ such that $f(u)\leq s$.
Now, $f(u)\leq q$ implies $f(u)\leq ^* q$ thus $F([u])=[f(u)]\trianglelefteq [q]$ and $F([u])\in F[ \CD]$.

Thus the partial order $\sq \P$ is quasi homogeneous indeed.
\kdok
\begin{cla}   \label{T6031}
The Boolean completion of a separative quasi homogeneous partial ordering is a homogeneous  complete Boolean algebra.
\end{cla}
\dok
Let $\P =\la P\leq\ra$ be a separative quasi homogeneous partial order.
First we show that
\begin{equation}\label{EQ6049}
\forall p\in P \;\;\ro \P \cong \ro (p\down ) .
\end{equation}
If $p\in P$, then there is a dense subset $D$ of $\P$ and an embedding $f:D\hookrightarrow p\down$ such that
$f[D]$ is a dense subset of $p\down$. Thus $D$ and $f[D]$ are isomorphic separative posets, which implies that $\ro D\cong \ro f[D] $.
In addition,  $D$ is a dense suborder of the separative order $\P$, which, by the uniqueness of the Boolean completion, implies
$\ro \P \cong \ro D$ and, similarly, $\ro f[D]\cong \ro (p\down )$
 and (\ref{EQ6049}) is true.

Let $\B =\ro \P$, $b\in \B ^+$ and w.l.o.g.\ suppose that $\P$ is a dense suborder of $\B ^+$.
Then there is $p\in \P$ such that $p\leq _\B b$.
Clearly the set
$(p\down)_\B \cap P = (p\down)_\P $ is a dense suborder of the relative algebra $(p\down)_\B$,
which implies
$(p\down)_\B \cong \ro ((p\down)_\P )$ so, by (\ref{EQ6049}),
$(p\down)_\B \cong \ro \P \cong \B$.
It is well known that, if $\B$ is a $\s$-complete Boolean algebra, $a,b\in \B$, $a\leq b$ and $\B\cong a\down$, then $\B\cong b\down$ (see \cite{Kopp1}, p.\ 180).
So we have $b\down \cong \B$.
\kdok
\begin{ex}\rm\label{EX6011}
Clearly homogeneous partial orders are quasi homogeneous, but the converse is not true. Let $\R$ be the real line and
$$
\P =\Big\la \{ (a,b] : a,b\in \R \land a<b \} \cup \{ \R \}, \subset \Big\ra .
$$
Then for $p=(a,b]$ we have $p\down \not\cong \P$, since the largest element of $\P$ is not the supremum of two smaller elements.
Thus $\P$ is not a homogeneous partial order.
On the other hand, if $f:\R \rightarrow (a, b)$ is an isomorphism, then it is easy to show that
the mapping $F:P\rightarrow p\down$ defined by $F(\R )=p$ and $F((c,d])=(f(c),f(d)]$ is an embedding and that
$F[P]$ is a dense subset of $p\down$. Thus the partial order $\P$ is quasi homogeneous. We note that $\P$ is, in addition, separative.
\end{ex}
\begin{te}\label{T4013}
For each relational structure $\X$ the Boolean completion $\ro \sq  \P (\X )$ of the poset $\P (\X )$ is a homogeneous complete Boolean algebra,
forcing equivalent to $\P (\X )$.
All generic extensions by $\P (\X )$ are elementarily equivalent.
\end{te}
\dok
By Theorem 2.2 of \cite{Ktow} the poset $\P (\X )$ is homogeneous and, by Theorem \ref{T6030}, $\ro \sq  \P (\X )$ is a homogeneousthe algebra.
By Fact \ref{T4042}(b) the posets $\P (\X )$ and $\ro \sq  \P (\X )$ are forcing equivalent.
By Theorem 4.3 of \cite{Ktow} either $|\sq \P (\X )|=1$, and then all generic extensions are trivial, or $\sq \P (\X )$ is an atomless poset, and then
$\B := \ro \sq  \P (\X )$ is an infinite homogeneous complete Boolean algebra.
This implies that for each $a,b\in \B \setminus \{ 0,1 \}$ there is $f\in \Aut (\B)$ such that $f(a)= b$
(see \cite{Kopp1}, Proposition 9.13) and, hence $\B ^+$ is a weakly homogeneous
partial order (we recall that a partial order $\P =\la P ,\leq \ra$ is called {\it weakly homogeneous} iff for each $p,q\in P$ there is $f\in \Aut (\P)$ such that $f(p)\not\perp q$).
By a known fact concerning weakly homogeneous partial
orders (see \cite{Kun}, p.\ 245), for each
sentence $\f$ of the language of set theory we have $1\Vdash \f$ or  $1\Vdash \neg \f$. Thus all generic extensions by $\P (\X )$
satisfy the same set of sentences.
\kdok
\subsection{Forcing-equivalence and isomorphism of Boolean completions}
Here we show that the posets of copies of two structures are forcing equivalent iff
their Boolean completions are isomorphic.
\begin{fac}   \label{T6002}
If $\B$ and $\C$ are complete Boolean algebras such that some $\B$-generic extension is equal to some $\C$-generic extension, then

(a) There are $b\in \B$ and $c\in \C$ such that $b \down \cong c\down$ (see \cite{Jech}, p.\ 267);

(b) If $\B$ and $\C$ are homogeneous algebras, then $\B \cong \C$. So, $\B \equiv \C \Leftrightarrow \B \cong \C$.
\end{fac}
\dok
(b) If $b\in \B$ and $c\in \C$ are the elements from (a), by the homogeneity we have $\B\cong b\down$ and $\C \cong c\down$ and, hence, $\B \cong \C$.
\kdok
\begin{te}\label{T6032}
Let $\X$ and $\Y$ be arbitrary relational structures. Then

(a) $\P (\X )\equiv \P (\Y )$ iff $\;\ro \sq \P (\X )\cong \ro \sq \P (\Y )$;

(b) The collections $\ge \P (\X )$ and $\ge \P (\Y )$ are either disjoint or equal.
\end{te}
\dok (a) By Fact \ref{T4042}(b), Fact \ref{T6002}(b) and Theorem \ref{T4013},
$\P (\X ) \equiv  \P (\Y )$ iff $\ro \sq \P (\X )$ $ \equiv \ro \sq \P (\Y )$ iff
$\ro \sq \P (\X ) \cong \ro \sq  \P (\Y )$.

(b) If $\ge \P (\X ) \cap \ge \P (\Y )\neq \emptyset$, then by Fact \ref{T4042}(b)
and Fact \ref{T6002}(b) we have $\ro \sq \P (\X )\cong \ro \sq \P (\Y )$, which implies
$\ge (\ro \sq \P (\X ))= \ge (\ro \sq \P (\Y ))$, that is $\ge  \P (\X )= \ge \P (\Y )$.
\kdok
\subsection{Isomorphic structures, equimorphic structures}
In this section we prove that the posets of copies of isomorphic (resp.\ equimorphic) structures are isomorphic
(resp.\ have isomorphic Boolean completions).
We will use the following elementary  fact.
\begin{fac}   \label{T3001}
Let $\la \P, \leq \ra$ be a pre-order and $p\in \P$. Then

(a) If $G$ is a $\P$-generic filter over $V$ and $p\in G$, then $G\cap p\down $ is a  $p\down$-generic filter over $V$
    and $V_\P [G]= V_{p\downarrow }[G\cap p\down]$;

(b) If $H$ is a $p\down$-generic filter over $V$, then $H\up$ is a $\P$-generic filter over $V$ and
    $V_{p\downarrow }[H]=V_\P [H\up ]$.
\end{fac}
\begin{lem}\label{T6005}
If $\X$ and $\Y$ are structures of the same language, $h: \X \hookrightarrow \Y$, and $C=h[X]$, then the mapping
$F: \P (\X ) \rightarrow (C\down )_{\P (\X , \Y )} $ defined by $F(A)=h[A]$, for $A\in \P (\X )$, is an isomorphism of the posets
$\la \P (\X ), \subset \ra$ and $\la (C\down )_{\P (\X , \Y )} ,\subset \ra $.
\end{lem}
\dok
For $A\in \P (\X)$ there is $\f :\X \hookrightarrow \X$ such that $\f [X]=A$ and, clearly,
$h\circ \f : \X \hookrightarrow \Y$, thus $h[\f [X]]=h[A]\in \P (\X , \Y )$ and $h[A]\subset h[X]=C$, which implies that $h[A]\in (C\down )_{\P (\X , \Y )}$.
So $F[\P (\X )] \subset (C\down )_{\P (\X , \Y )} $.

Since $h$ is an injection, for each $A,B\in \P (\X )$ we have $F(A)\subset F(B)$ iff $h[A]\subset h[B]$ iff
$h^{-1}[h[A]]\subset h^{-1}[h[B]]$ iff $A\subset B$, thus $F$ is an embedding.

If $D\in \P (\X , \Y )$ and $D\subset C$, then $h[h^{-1}[D]]=D$ and
the surjective restriction $h\mid h^{-1}[D] : h^{-1}[D]\rightarrow D$ is an isomorphism, which implies $h^{-1}[D]\in \P (\X )$. In addition $F(h^{-1}[D])=h[h^{-1}[D]]=D$ thus $F$ is onto.
\kdok
\begin{te}\label{T6004}
If $\X$ and $\Y$ are structures of the same relational language, then

(a) $\X\cong \Y  \;\Rightarrow \; \P (\X ) \cong \P (\Y )$;

%(b) $\X \rightleftarrows \Y \;\Rightarrow \;\P (\X )\equiv \P (\Y)$.

(b) $\X \rightleftarrows \Y \;\Rightarrow \;\ro\sq \P (\X )\cong \ro\sq \P (\Y)$.
\end{te}
\dok
(a) If $h:\X \rightarrow \Y$ is an isomorphism, then, by Lemma \ref{T6005},
$\la \P (\X ), \subset \ra \cong \la \P (\X ,\Y ), \subset \ra $
and, clearly, $\la \P (\X ,\Y ), \subset \ra =\la \P (\Y ), \subset \ra$.

(b) Let $f:\X \hookrightarrow \Y$, $g:\Y \hookrightarrow \X$ and $\P (\Y )\up =\{ S\subset Y : \exists B\in \P (\Y )\; B\subset S \}$.

First we show that $\P (\X , \Y ) := \{ h[X] \mid h: \X \hookrightarrow \Y \} $ is a dense suborder of $\la \P (\Y )\up , \subset \ra$.
If $C\in \P (\X , \Y )$ and $h: \X \hookrightarrow \Y$, where $C=h[X]$, then, clearly,
$h\circ g : \Y \hookrightarrow \Y$  and, hence, $h[g[Y]]\in \P (\Y )$ and $h[g[Y]]\subset h[X]=C$, which implies $C\in \P (\Y )\up$. Thus
$\P (\X , \Y )\subset \P (\Y )\up$. Let $S\in \P (\Y )\up$, $B\in \P (\Y)$, where $B\subset S$ and $\p : \Y \hookrightarrow \Y$, where $B=\p [Y]$.
Now $\p \circ f :\X \hookrightarrow \Y$ and, hence, $\p[ f [X ]]\in \P (\X , \Y )$ and $\p[ f [X ]] \subset \p [Y]=B\subset S$. Thus
$\P (\X , \Y )$ is dense in $\la \P (\Y )\up , \subset \ra$. Since $\P (\Y )$ is dense in $\la \P (\Y )\up , \subset \ra$ as well, we have
\begin{equation}\label{EQ6002}
\la \P (\X , \Y ) ,\subset \ra \equiv \la \P (\Y ), \subset \ra .
\end{equation}
Now let $W$ be a generic extension by $\P (\Y )$. By (\ref{EQ6002}) $W=V_{\P (\X , \Y )}[G]$, where $G$ is a $\P (\X , \Y )$-generic filter over $V$.
Let $C\in G$. By Fact \ref{T3001}(a) we have $V_{\P (\X , \Y )}[G]= V_{C\down }[G\cap C\down ]$ and,
if $F:\la \P (\X ), \subset \ra \rightarrow \la (C\down )_{\P (\X , \Y )} ,\subset \ra $ is the isomorphism defined in Lemma \ref{T6005},
then $H:= F^{-1}[G\cap C\down]$ is a $\P (\X )$-generic filter over $V$ and $V_{C\down }[G\cap C\down ]=V_{\P (\X )}[H]$. Thus $W=V_{\P (\X )}[H]$ and, by
Theorem \ref{T6032}(b), $\P (\X ) \equiv \P (\Y )$. Now, by Theorem \ref{T6032}(a),  $\ro\sq \P (\X )\cong \ro\sq \P (\Y)$.
\kdok
\section{The hierarchy of similarities on the class Mod$_L$(X)}\label{S3}
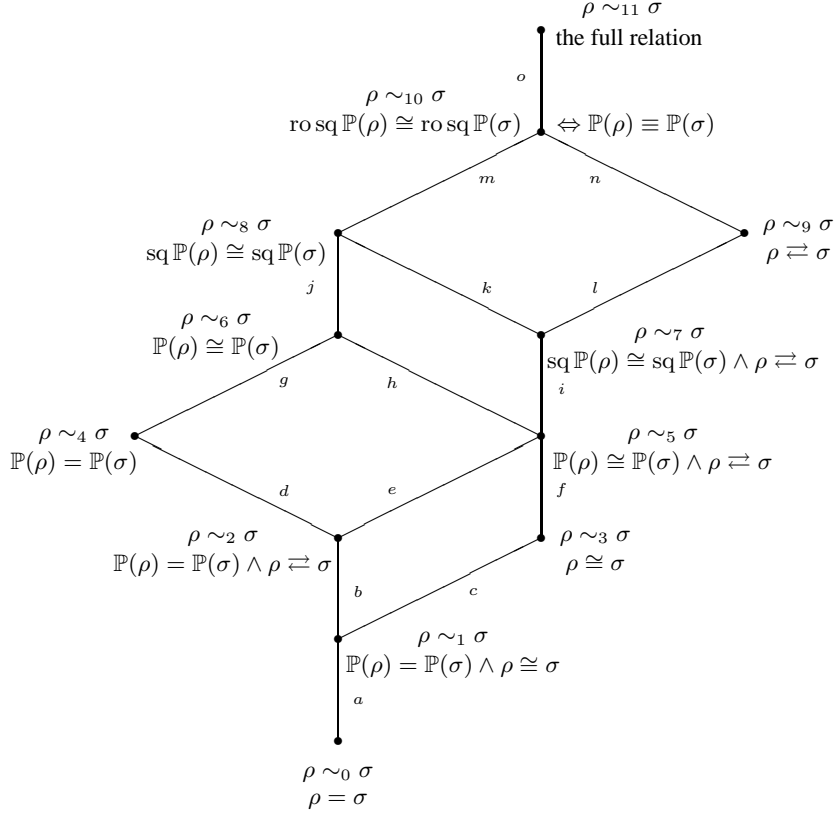
\begin{figure}[htb]%\label{F4008}
\unitlength 0.9mm %1mm % = .854pt
\linethickness{0.5pt}
\ifx\plotpoint\undefined\newsavebox{\plotpoint}\fi % GNUPLOT compatibility
\begin{picture}(140,122)(0,0)

%----------------------------- linije --------------------------

\put(55,10){\line(0,1){15}}%1
\put(55,25){\line(0,1){15}}%2
\put(55,25){\line(2,1){30}}%3
\put(55,40){\line(-2,1){30}}%4
\put(55,40){\line(2,1){30}}%5
\put(85,40){\line(0,1){15}}%6
\put(25,55){\line(2,1){30}}%7
\put(85,55){\line(-2,1){30}}%8
\put(85,55){\line(0,1){15}}%9
\put(55,70){\line(0,1){15}}%10
\put(85,70){\line(-2,1){30}}%11
\put(85,70){\line(2,1){30}}%12
\put(55,85){\line(2,1){30}}%13
\put(115,85){\line(-2,1){30}}%14
\put(85,100){\line(0,1){15}}%15

%----------------------------- tacke --------------------------

\put(55,10){\circle*{1}}%1
\put(55,25){\circle*{1}}%
\put(55,40){\circle*{1}}%
\put(85,40){\circle*{1}}%
\put(25,55){\circle*{1}}%
\put(85,55){\circle*{1}}%
\put(55,70){\circle*{1}}%
\put(85,70){\circle*{1}}%
\put(55,85){\circle*{1}}%
\put(115,85){\circle*{1}}%
\put(85,100){\circle*{1}}%
\put(85,115){\circle*{1}}%

%----------------------------- tekst -----------------------
%\small
\footnotesize
%\scriptsize
%\tiny

\put(55,5){\makebox(0,0)[cc]{$\r \sim _0 \s $}}%0
\put(55,1){\makebox(0,0)[cc]{$  \r=\s$}}%0'
\put(72,25){\makebox(0,0)[cc]{$\r \sim _1 \s $}}%1
\put(72,21){\makebox(0,0)[cc]{$\P (\r )=\P (\s ) \land \r \cong \s$}}%1'
\put(38,40){\makebox(0,0)[cc]{$\r \sim _2 \s $}}%2
\put(38,36){\makebox(0,0)[cc]{$ \P (\r )=\P (\s ) \land \r \rightleftarrows \s$}}%2'
\put(93,40){\makebox(0,0)[cc]{$\r \sim _3 \s$}}%3
\put(93,36){\makebox(0,0)[cc]{$\r \cong \s $}}%3'
\put(16,55){\makebox(0,0)[cc]{$\r \sim _4 \s $}}%4
\put(16,51){\makebox(0,0)[cc]{$ \P ( \r )=\P (\s )$}}%4'
\put(103,55){\makebox(0,0)[cc]{$\r \sim _5 \s $}}%5
\put(103,51){\makebox(0,0)[cc]{$\P (\r )\cong \P (\s ) \land \r \rightleftarrows \s$}}%5'
\put(37,72){\makebox(0,0)[cc]{$\r \sim _6 \s $}}%6
\put(37,68){\makebox(0,0)[cc]{$ \P ( \r )\cong \P (\s )$}}%6'
\put(104,70){\makebox(0,0)[cc]{$\r \sim _7 \s $}}%7
\put(106,66){\makebox(0,0)[cc]{$\sq \P (\r )\cong \sq \P (\s )\land \r \rightleftarrows \s$}}%7'
\put(40,86){\makebox(0,0)[cc]{$\r \sim _8 \s$}}%8
\put(40,82){\makebox(0,0)[cc]{$\sq \P ( \r )\cong \sq \P (\s )$}}%8'
\put(123,86){\makebox(0,0)[cc]{$\r \sim _9 \s$}}%9
\put(123,82){\makebox(0,0)[cc]{$\r \rightleftarrows \s$}}%9'
\put(65,105){\makebox(0,0)[cc]{$\r \sim _{10} \s$}}%10
\put(65,101){\makebox(0,0)[cc]{$\ro \sq  \P ( \r )\cong \ro  \sq \P (\s )$}}%10'
\put(99,101){\makebox(0,0)[cc]{$\Leftrightarrow\P ( \r )\equiv \P (\s )$}}%10'
\put(98,114){\makebox(0,0)[cc]{the full relation}}%A
\put(97,118){\makebox(0,0)[cc]{$\r \sim _{11} \s$}}%10
{\tiny
\put(58,16){\makebox(0,0)[cc]{$a$}}%
\put(58,32){\makebox(0,0)[cc]{$b$}}%
\put(75,32){\makebox(0,0)[cc]{$c$}}%
\put(47,47){\makebox(0,0)[cc]{$d$}}%
\put(63,47){\makebox(0,0)[cc]{$e$}}%
\put(88,47){\makebox(0,0)[cc]{$f$}}%
\put(47,63){\makebox(0,0)[cc]{$g$}}%
\put(63,63){\makebox(0,0)[cc]{$h$}}%
\put(88,62){\makebox(0,0)[cc]{$i$}}%
\put(51,77){\makebox(0,0)[cc]{$j$}}%
\put(77,77){\makebox(0,0)[cc]{$k$}}%
\put(93,77){\makebox(0,0)[cc]{$l$}}%
\put(77,93){\makebox(0,0)[cc]{$m$}}%
\put(93,93){\makebox(0,0)[cc]{$n$}}%
\put(82,108){\makebox(0,0)[cc]{$o$}}%
}
\end{picture}
\caption{Some implications between the similarities on $\Mod _L (X)$}\label{F4008}
\end{figure}

Now we restrict our consideration to some smaller classes of structures. If $L=\la R_i : i\in I\ra$ is a language, $X$ a fixed set and
$\r=\la \r _i :i\in I \ra\in \Int _L (X)$, we will abuse notation writing $\P (\r )$ instead of $\P (\la X, \r\ra)$ and $\la \P (\la X, \r\ra),\subset \ra$
whenever the context admits it.
So, restricting our similarity relations to the set $\Mod _L (X)$ or, equivalently,
to the corresponding set of interpretations, $\Int _L (X)$, we obtain the following equivalence relations:
for $\r=\la \r _i :i\in I \ra, \s =\la \s _i :i\in I \ra\in \Int _L (X)$ (writing $\r \cong \s$ instead of $\la X, \r \ra \cong\la X, \s \ra$
and similarly for $\r \rightleftarrows \s$) we define
\begin{center}
{\tabcolsep 8pt
\begin{tabular}{ll}
$\r \sim _0 \s \Leftrightarrow  \r=\s$
&
$\r \sim _6 \s \Leftrightarrow \P ( \r )\cong \P (\s )$
\\
$\r \sim _1 \s \Leftrightarrow \P ( \r )=\P (\s ) \land \r \cong \s$
&
$\r \sim _7 \s \Leftrightarrow \sq \P ( \r )\cong \sq \P ( \s )\land \r \rightleftarrows \s$
\\
$\r \sim _2 \s \Leftrightarrow \P ( \r )=\P (\s ) \land \r \rightleftarrows \s$
&
$\r \sim _8 \s \Leftrightarrow \sq \P ( \r )\cong \sq \P ( \s )$
\\
$\r \sim _3 \s \Leftrightarrow \r \cong \s$
&
$\r \sim _9 \s \Leftrightarrow \r \rightleftarrows \s$
\\
$\r \sim _4 \s \Leftrightarrow \P ( \r )=\P (\s )$
&
$\r \sim _{10} \s \Leftrightarrow \P ( \r )\equiv \P ( \s )$
\\
$\r \sim _5 \s \Leftrightarrow \P ( \r )\cong \P (\s ) \land \r \rightleftarrows \s$
&
$\r \sim _{11} \s \Leftrightarrow 0=0$.
\end{tabular}
}
\end{center}
Then some implications between the similarities on the set $\Mod _L (X)$ are displayed in Figure \ref{F4008}.

\begin{figure}[htb]%\label{F4010}
\begin{center}
\unitlength 1mm %0.7mm % = .854pt
\linethickness{0.5pt}
\ifx\plotpoint\undefined\newsavebox{\plotpoint}\fi % GNUPLOT compatibility
\begin{picture}(80,30)(0,0)
%----------------------------- linije --------------------------
\put(10,5){\line(0,1){20}}%1
%----------------------------- tacke --------------------------
\put(10,5){\circle*{1}}%1
\put(10,15){\circle*{1}}%
\put(10,25){\circle*{1}}%
%----------------------------- tekst -----------------------
\small
%\footnotesize
%\scriptsize
%\tiny
\put(48,25){\makebox(0,0)[cc]{$\sim _4 \;=\; \sim _6 \;=\; \sim _8 \;=\; \sim _{10} \;=\;\sim _{11} \;=\;$ the full relation}}%
\put(52,15){\makebox(0,0)[cc]{$\sim _1 \;=\; \sim _2 \;=\; \sim _3 \;=\; \sim _5 \;=\; \sim _7 \;=\; \sim _9 \;=\;$ the isomorphism}}%
\put(26,5){\makebox(0,0)[cc]{$\sim _0 \;=\;$ the equality}}%1
\end{picture}
\end{center}
\caption{The similarities on the class $\Mod _L (X)$, if $1<|X|<\o$}\label{F4010}
\end{figure}
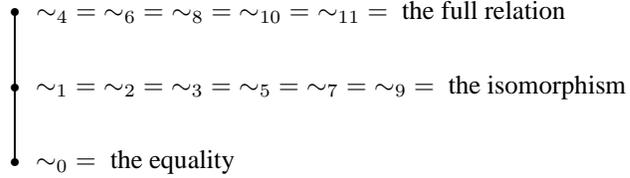

It is natural to ask
are there more implications in it (except the ones which follow from the transitivity), that is,
are some of the implications $a$ - $o$, in fact,  equivalences.
Concerning this question we will show that the class of all relational structures splits into the following three parts:
finite structures, infinite structures of unary languages, and infinite structures of non-unary languages.
(A language $L=\la R_i : i\in I\ra$ is  called {\it unary} iff $\ar (R_i)=1$, for all $i\in I$. Structures of unary languages
will be called {\it unary structures}).
Let us call a class $\CC$ of structures a {\it Cantor-Schr\"{o}der-Bernstein (CSB) class} iff
$$
\forall \X ,\Y \in \CC \;\; (\X \rightleftarrows \Y \;\Rightarrow \; \X \cong \Y ).
$$
For finite structures the diagram from Figure \ref{F4008} collapses significantly.
\begin{ex}\rm
If $L$ is an arbitrary relational language and $X$ a finite set, then for each $\r \in \Int _L(X)$  we have $\P (\r )=\{ X \}$, because
$X\in \P (X,\r )\subset [X]^{|X|}=\{ X \}$. Thus, $\sim_4$ is the full relation, which implies that $\sim _4 = \sim _6 = \sim _8 = \sim _{10} = \sim _{11}$.
In addition, $\Mod _L (X)$ is a CSB class. Namely, if $\r ,\s \in \Int _L (X)$ and $\r  \rightleftarrows \s$, then there is an embedding
$f:\la X, \r \ra \rightarrow \la X, \s \ra$, and, since $X$ is a finite set, $f$ is an isomorphism, thus $\r \cong \s$.
So we have $\sim _9 \subset \sim _1 $, which implies
$\sim _1 = \sim _2 = \sim _3 = \sim _5 = \sim _7 = \sim _9$.
Since $\la X, \la \emptyset , \emptyset , \dots \ra\ra \not\cong \la X, \la X^{n_i}: i\in I \ra\ra$, we have
$\sim _3 \neq \sim _{11}$. If $|X|>1$, let $a$ and $b$ be different elements of $X$, $i_0\in I$,
and let $\r , \s \in \Int _L (X)$, where
$\r _{i_0}= \{ \la a,a,\dots ,a\ra\} , \s _{i_0}= \{ \la b,b,\dots ,b\ra\}\subset X^{n_{i_0}}$ and
$\r_i =\s _i =\emptyset$, for $i\neq i_0$. Then $\r \neq _0 \s$, but $\r \cong \s$ and, hence, $\sim _0 \neq\sim _1$.
Thus Figure \ref{F4010} describes the hierarchy of the similarities $\sim _k$ on the set $\Mod _L (X)$, if $|X|>1$.

We prove that $\sim _0 =\sim _1 \;\Leftrightarrow \;|X|=1$. Let $X=\{ x \}$ and $\r , \s \in \Int _L (X)$, where $\r \sim _1 \s$.
Then there is an isomorphism $f: \la \{ x \} ,\r \ra \rightarrow\la \{ x \} ,\s \ra$ and, consequently, for each $i\in I$ we have
$\la x,x,\dots ,x\ra \in \r _i \Leftrightarrow \la x,x,\dots ,x\ra \in \s _i$ and, hence, $\r_i =\s _i$. So  $\r = \s$, that is $\r \sim _0 \s$ and the
inclusion $\sim _1 \subset \sim _0$ is proved.
\end{ex}

\subsection{Infinite unary structures}
In this subsection we assume that $L=\la R_i :i\in I\ra$ is a unary relational language.  If $\X =\la X, \la \r _i : i\in I \ra \ra$ is an $L$-structure,
it is easy to check that the binary relation $\approx$ on the set $X$ defined by:
$x\approx y \Leftrightarrow \forall i\in I \;\;(x\in \r_i \Leftrightarrow y\in \r_i )$
is an equivalence relation.  Then $[x]:=\{ y\in X : y\approx x\}$ is the equivalence class of $x\in X$, and if $X/\!\approx \;=\{ X_j :j\in J\}$
is the corresponding partition we define $\k _j:=|X_j|$, for $j\in J$, and $J_0 :=\{ j\in J: |X_j|<\o \}$.
\begin{te}\label{T6020}
Let $\X =\la X, \la \r _i : i\in I \ra \ra$ be a unary structure. Then

(a) If $f:X\rightarrow X$ is an injection, then $f\in \Emb (\X ) \Leftrightarrow \forall x\in X \;\; f[[x]] \subset [x]$;

(b) If $J_0 =J$, then $\P (\X )=\{ X \}$;

(c) If $J_0 \neq J$,  then the poset $\P (\X )$ is atomless and we have
$$\textstyle
\P (\X ) \cong \prod _{j\in J\setminus J_0}\la[\k _j]^{\k _j}, \subset \ra \;\;\mbox{ and }\;\;
\sq  \P (\X ) \cong \prod _{j\in J\setminus J_0}(P(\k _j )/[\k _j]^{<\k _j })^+ .
$$
\end{te}
\dok
(a) If $f \in \Emb (\X )$ and $x\in X$, then $x\in \r _i \Leftrightarrow f(x)\in \r _i$, for each
$i\in I$, thus $x\approx f(x)$. So for $y\in [x]$ we have $f(y)\approx y\approx x$ and, hence, $f(y)\in [x]$.

Let $f[[x]] \subset [x] $, for all $x\in X$. Then, for $x\in X$,  $x\in [x]$ implies $f(x)\in [x]$, that is $f(x)\approx x$. Hence
$\forall i\in I \; \forall x\in X \;(x\in \r _i \Leftrightarrow f(x)\in \r _i)$, so $f \in \Emb (\X )$.

(b) Let $J_0 =J$. By (a), for $f \in \Emb (\X )$ and $x\in X$ we have $f[[x]]\subset [x]$ and, since
$|[x]|< \o$, $f[[x]]=[x]$,  which implies $f[X]=X$.

\begin{figure}[htb]
\begin{center}
\unitlength 1mm % = .854pt
\linethickness{0.5pt}
\ifx\plotpoint\undefined\newsavebox{\plotpoint}\fi % GNUPLOT compatibility
\begin{picture}(80,50)(0,0)

%----------------------------- linije --------------------------

\put(30,5){\line(0,1){10}}%1
\put(30,15){\line(-1,1){10}}%2
\put(30,15){\line(1,1){10}}%3
\put(20,25){\line(1,1){10}}%4
\put(40,25){\line(-1,1){10}}%5
\put(30,35){\line(0,1){10}}%6
%----------------------------- tacke --------------------------

\put(30,5){\circle*{1}}%1
\put(30,15){\circle*{1}}%
\put(20,25){\circle*{1}}%
\put(40,25){\circle*{1}}%
\put(30,35){\circle*{1}}%
\put(30,45){\circle*{1}}%

%----------------------------- tekst -----------------------

\footnotesize
%\scriptsize
%\tiny
\put(48,45){\makebox(0,0)[cc]{$\sim _{11} \;=\; $ the full relation}}%
\put(38,15){\makebox(0,0)[cc]{$\sim _1 \;=\; \sim _2 $}}%
\put(16,26){\makebox(0,0)[cc]{$\sim _4 $}}%
\put(4,22){\makebox(0,0)[cc]{$=\;$ the equality of $\P (\X) $}}%
\put(54,35){\makebox(0,0)[cc]{$\sim _6 \; =\;$ the isomorphism of $\P (\X) $}}%
\put(57,26){\makebox(0,0)[cc]{$\sim _3 \;=\; \sim _5 \;=\; \sim _7 \;=\; \sim _9 $}}%
\put(68,22){\makebox(0,0)[cc]{$\;=\;$ the isomorphism $\;=\;$  the equimorphism}}%
\put(44,5){\makebox(0,0)[cc]{$\sim _0 \; =\; $ the equality}}%1
\end{picture}
\end{center}
\caption{The similarities on $\Mod _L(\k )$, for unary $L$ and infinite $\k$}\label{F4011}
\end{figure}
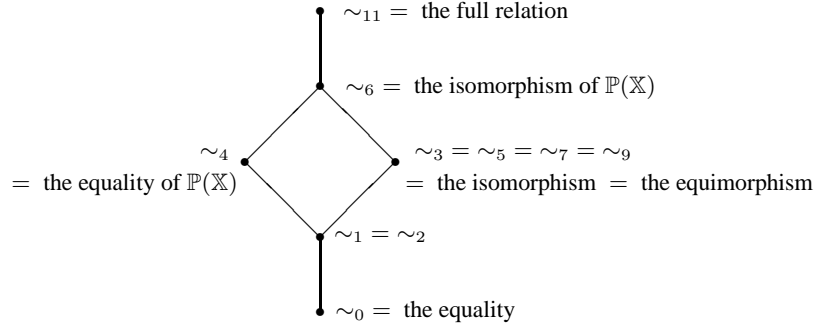

(c) If $f \in \Emb (\X )$, then, by (a), $f[X_j]=X_j$, for all $j\in J_0$, and $C_j:= f[X_j]\in [X_j]^{\k _j}$, for all $j\in J\setminus J_0$. Thus
the inclusion ``$\subset$" in the equality
\begin{equation}\label{EQ6028}\textstyle
\P (\X )=\Big\{ \bigcup _{j\in J_0}X_j \cup \bigcup _{j\in J\setminus J_0}C_j : \la C_j :j\in J\setminus J_0\ra \in \prod _{j\in J\setminus J_0}[X_j]^{\k _j}\Big\}
\end{equation}
is proved. On the other hand, if $\la C_j :j\in J\setminus J_0\ra \in \prod _{j\in J\setminus J_0}[X_j]^{\k _j}$ and if we choose bijections
$\varphi _j : X_j \rightarrow C_j$, for all $j\in J\setminus J_0$, then by (a) we have
$f=\bigcup _{j\in J_0}\id _{X_j} \cup \bigcup _{j\in J\setminus J_0}\varphi _j \in \Emb (\X )$
and, hence, $\bigcup _{j\in J_0}X_j \cup \bigcup _{j\in J\setminus J_0}C_j\in \P (\X )$, so (\ref{EQ6028}) is true.
Thus the mapping $F:\prod _{j\in J\setminus J_0}\la [X_j]^{\k _j}, \subset \ra \rightarrow \la \P (\X ),\subset \ra$ given by
$$\textstyle
F(\la C_j :j\in J\setminus J_0\ra )=\bigcup _{j\in J_0}X_j \cup \bigcup _{j\in J\setminus J_0}C_j
$$
is a well-defined surjection and, since $\{ X_j : j\in J \}$ is a partition of $X$, it is an injection. It is easy to see
that $F$ is an order isomorphism. By Fact \ref{T4042}(c) we have
$\sq \la \P (\X ), \subset \ra \cong \prod _{j\in J\setminus J_0}\sq \la[\k _j]^{\k _j}, \subset \ra =
\prod _{j\in J\setminus J_0}(P(\k _j )/[\k _j]^{<\k _j })^+ $.
\kdok
\begin{lem}\label{T6023}
Let $\k\geq \o$ be a cardinal, $U\subset \k$ and $\l :=\min \{ |U|,|\k \setminus U |\}$. Then
$\r =\la U, \emptyset, \emptyset, \dots \ra\in \Int _L (\k )$ and we have

(a) $\P ( \r ) =
\{ C_1 \cup C_2 : C_1 \in [U] ^{|U |} \land C_2 \in [\k \setminus U ]^{|\k \setminus U|}\} $;

(b) $\P ( \r )\cong
\la [\k ] ^\k ,\subset \ra \times \la [\l ]^\l ,\subset \ra $;

(c) $\sq  \P ( \r ) \cong
(P(\k )/[\k]^{<\k })^+ \times (P(\l )/[\l ]^{<\l })^+$, where,  by convention,
for $\l \in \o$, by $(P(\l )/[\l]^{<\l })^+$ we denote the one-element poset.
\end{lem}
\dok
For $x,y\in \k$ we have: $x\approx y$ iff $x\in \r _i \Leftrightarrow y\in \r _i$, for all $i\in I$,
iff $x\in U \Leftrightarrow y\in U$. Thus $\k /\approx =\{ U, \k \setminus U \} $ and we apply Theorem \ref{T6020}.
\kdok
\begin{fac}\label{T6025}
(a) If $\k>\o$ is a regular cardinal and $2^\k =\k ^+$, then $\ro (P(\k )/[\k ]^{<\k })\cong \Col (\o , 2^\k)$ (Balcar, Vop\v enka \cite{BalcVop}; see also \cite{BalcSim}, p.\ 380).

(b) Under CH, all separative atomless $\o _1$-closed posets of size $\o _1$ are forcing equivalent (for example to $\Col (\o _1 ,\o _1 )$) (folklore).

(c)  If $\l >\o$ is a cardinal and $\P$ a poset of size $\l$ such that $1_\P \Vdash |\check{\l }|=\check{\o } $,
then $\ro \P \cong \Col (\o ,\l )$ (see \cite{Jech}, p.\ 277).

(d) If $\B$ is a Boolean algebra of size $>2$, then $\B ^+ \not\cong \B ^+\times \B ^+$.
\end{fac}
\dok
(d) The sentence $\forall x\neq 1 \; \exists _1 y \; (x\perp y \;\&\; x\vee y =1)$ is true in the poset $\B ^+ ,$ but it is not true
in  $\B ^+\times \B ^+$. Namely, since $|\B |>2$, there is $a\in \B ^+ \setminus \{ 1 \}$ and we have
$x:= \la 1, a \ra \in (\B ^+\times \B ^+ )\setminus \{\la 1,1 \ra\}$ and $a'\in \B ^+$, but for each $b\in \B ^+$ we have
$\la 1,a \ra \perp \la b, a'\ra$ and $\la 1,a \ra \vee \la b, a'\ra =\la 1,1 \ra$.
\kdok
\begin{te} \label{T6024}
For any unary language $L$ and infinite cardinal $\k$ we have

(a) $\Mod _L (\k )$ is a CSB class;

(b) Figure \ref{F4011} describes the hierarchy of the similarities $\sim _k$, for $k\neq 8,10$, on the set $\Mod _L (\k )$.
In addition we have $\sim _8 \;\neq \;\sim _{11}$.

(c) If $\k$ is a regular cardinal and $2^\k =\k ^+$, then $\sim _8 \;\neq \;\sim _{10}$.
\end{te}

\dok
Let  $L=\{ R_i: i\in I \}$.

(a)  Assuming that $\r =\la \r _i :i\in I \ra , \s =\la \s _i :i\in I \ra \in \Int _L (\k )$ and $\r \rightleftarrows \s$ we show that
$\r \cong \s$. By the assumption, there are embeddings
\begin{equation}\label{EQ6032}\textstyle
f: \la \k ,\r \ra \hookrightarrow \la \k ,\s \ra \;\;\mbox{ and }\;\; g: \la \k ,\s \ra \hookrightarrow \la \k ,\r \ra.
\end{equation}
Let $\approx _\r$ and $\approx _\s$ be the equivalence relations  determined by the interpretations $\r$
and $\s$ respectively (see Theorem \ref{T6020}) and, for $x\in \k$, let $[x]_\r$ and $[x]_\s$ be the corresponding equivalence classes. First we prove that
\begin{equation}\label{EQ6033}\textstyle
\forall x\in \k \;\; f[[x]_\r ]\subset [f(x)]_\s \;\;\mbox{ and }\;\;\forall x\in \k \;\; g[[x]_\s ]\subset [g(x)]_\r .
\end{equation}
For a proof of the first statement we take $x\in \k$ and $y\in [x ]_\r$. Then $y\approx _\r x$, that is
\begin{equation}\label{EQ6035}\textstyle
\forall i\in I \;\; (x \in \r _i \Leftrightarrow y\in \r _i ) ,
\end{equation}
and, since $f$ is an embedding, we have
\begin{equation}\label{EQ6034}\textstyle
\forall i\in I \;\; \forall x\in \k\;\;(x\in \r _i \Leftrightarrow f(x)\in \s _i ) .
\end{equation}
We prove that $f(y)\in [f(x)]_\s$, which means that $f(y)\in \s _i \Leftrightarrow f(x)\in \s _i$, for all $i\in I$.
So, $f(y)\in \s _i$ iff (by (\ref{EQ6034})) $y\in \r _i$ iff (by (\ref{EQ6035})) $x\in \r_i$ iff (by (\ref{EQ6034})) $f(x)\in \s _i$.
Thus the first statement of (\ref{EQ6033}) is proved and the second has a symmetric proof.

Let $\k = \bigcup _{j\in J}X_j$ and $\k = \bigcup _{k\in K}Y_k$ be the partitions determined by the relations $\approx _\r$ and $\approx _\s$
respectively. By (\ref{EQ6033}), if $j\in J$ and $X_j =[x ]_\r$, then $f[X_j]\subset [f(x)]_\s =Y_k$, for (a unique) $k\in K$. Similarly,
for each $k\in K$ there is a unique $j\in J$ satisfying $g[Y_k]\subset X_j$ so we define the functions
\begin{equation}\label{EQ6036}\textstyle
F:J\rightarrow K \;\;\mbox{ by: }\;\; F(j)=k \;\mbox{ iff }\; f[X_j]\subset Y_k  ,
\end{equation}
\begin{equation}\label{EQ6037}\textstyle
G:K\rightarrow J \;\;\mbox{ by: }\;\; G(k)=j \;\mbox{ iff }\; g[Y_k]\subset X_j ,
\end{equation}
and prove that
\begin{equation}\label{EQ6038}\textstyle
G\circ F =\id _J \;\;\mbox{ and }\;\; F\circ G =\id _K.
\end{equation}
By (\ref{EQ6032}) we have $g\circ f:\la \k ,\r \ra \hookrightarrow \la \k ,\r \ra$ and, by Theorem \ref{T6020}(a),
\begin{equation}\label{EQ6039}\textstyle
\forall x\in X \;\; g\big[f\big[[x]_\r \big]\big]\subset [x]_\r .
\end{equation}
For $j\in J$ we prove that $G(F(j))=j$. Let $F(j)=k$ and $x\in X_j$. Then $X_j=[x]_\r$, by (\ref{EQ6033})
$f[X_j]=f[[x]_\r ]\subset [f(x)]_\s =Y_{k'}$, for some $k'\in K$, and, by (\ref{EQ6036}) $f[X_j]\subset Y_k$, which implies $k'=k$. Thus $f[[x]_\r ]\subset Y_k$
and, hence,
\begin{equation}\label{EQ6040}\textstyle
g[f[[x]_\r ]]\subset g[Y_k ] .
\end{equation}
Let $G(k)=j'$. Then by (\ref{EQ6037}) and (\ref{EQ6040}) we have $g[f[[x]_\r ]]\subset g[Y_k ]\subset X_{j'}$ and, by (\ref{EQ6039}),
$g[f[[x]_\r ]]\subset X_j$, which implies $j'=j$. Thus $G(F(j))=G(k)=j$ and the first equality in (\ref{EQ6038}) is proved. The second equality has a similar proof.

Now we prove that
\begin{equation}\label{EQ6041}\textstyle
\forall j\in J \;\; |X_j|=|Y_{F(j)}| .
\end{equation}
By (\ref{EQ6036}) we have $|X_j|=|f[X_j]|\leq| Y_{F(j)}|$ and, by (\ref{EQ6037}) and (\ref{EQ6038}),
$|Y_{F(j)}|=|g[Y_{F(j)}]| \leq |X_{G(F(j))}| =|X_j|$. So (\ref{EQ6041}) is true.

By (\ref{EQ6041}) there are bijections $\f _j:X_j \rightarrow Y_{F(j)}$; let
$\f =\bigcup _{j\in J}\f _j : \k \rightarrow \k$. Since $\{ X_j : j\in J \}$ is a partition of $\k$ the mapping $\f$ is well defined.
By  (\ref{EQ6038}) $F:J\rightarrow K$ is a bijection and, since the mappings $\f _j$ are surjections, $\f$ is a surjection as well.
Since $\{ Y_k : k\in K \}$ is a partition of $\k$ and the mappings $\f _j$ are injections, $\f$ is a injection too.
Thus $\f$ is a bijection from $\k$ onto $\k$.

In order to show that $\f :\la \k ,\r \ra \rightarrow \la \k ,\s \ra$ is an isomorphism, that is
\begin{equation}\label{EQ6042}\textstyle
\forall i\in I \;\; \forall x\in \k\;\;(x\in \r _i \Leftrightarrow \f (x)\in \s _i ) ,
\end{equation}
we take $i_0 \in I$ and $x_0 \in \k$. Let $j\in J$, where $x_0\in X_j$. Then $X_j=[x_0]_\r$ and $\f (x_0)=\f _j (x_0 )\in Y_{F(j)}$ and, by
(\ref{EQ6033}) and (\ref{EQ6036}), $f(x_0)\in [f(x_0)]_\s =Y_{F(j)}$. Thus $\f (x_0 )\approx _\s f(x_0)$, that is
\begin{equation}\label{EQ6043}\textstyle
\forall i\in I \;\; (\f (x_0) \in \s _i \Leftrightarrow f(x_0)\in \s _i ) .
\end{equation}
Now $x_0 \in \r _{i_0 }$ iff (by (\ref{EQ6034})) $f(x_0) \in \s _{i_0 }$ iff (by (\ref{EQ6043})) $\f (x_0) \in \s _{i_0 }$ and (\ref{EQ6042}) is proved.
Thus $\f :\la \k ,\r \ra \rightarrow \la \k ,\s \ra$ is an isomorphism and, hence, $\r \cong \s$.

(b) By (a) we have $\sim _9 \;\subset \;\sim _3$ which, according to Figure \ref{F4008}, implies that $\sim _3 \;=\; \sim _5 \;=\; \sim _7 \;=\; \sim _9$ and
$\sim _1 \;=\; \sim _2$.

Let us prove that $\sim _0 \;\varsubsetneq \;\sim _1$.
If $\k =A\cup B$, where $A\cap B=\emptyset$ and $A,B\in [\k ]^\k$, then $\r :=\la A ,\emptyset, \emptyset ,\dots  \ra \neq
\s :=\la B ,\emptyset, \emptyset ,\dots  \ra$. By Lemma \ref{T6023}(a) we have
$\P ( \r )= \{ C_1 \cup C_2 : C_1 \in [A]^\k \land C_2 \in [B]^\k \} =\P ( \s )$. If $f:\k \rightarrow \k$ is a bijection
satisfying $f[A]=B$, then $f :\la \k ,\r \ra \rightarrow \la \k ,\s \ra$ is an isomorphism and, hence, $\r \sim _1 \s$, but $\r \not\sim _0 \s$.

Now we prove that $\sim _3 \;\not\subset \;\sim _4$ and, hence, $\sim _1 \;\varsubsetneq \;\sim _3$ and $\sim _4 \;\varsubsetneq \;\sim _6$.
Let $x,y\in \k$, $x\neq y$ and let $\r :=\la \{ x\},\emptyset, \emptyset ,\dots  \ra $ and $\s :=\la \{ y \} ,\emptyset, \emptyset ,\dots  \ra$.
If $f:\k \rightarrow \k$ is a bijection satisfying $f(x)=y$,
then $f :\la \k ,\r \ra \rightarrow \la \k ,\s \ra$ is an isomorphism and, hence, $\r \sim _3 \s$.
By Lemma \ref{T6023}(a) we have
$$\P ( \r )= \{ C_1 \cup C_2 : C_1 \in [\{ x \}]^1 \land C_2 \in [\k \setminus \{ x \}]^\k \}=\{ C\in [\k ]^\k :x\in C\}$$
and, similarly,
$\P (\s )=\{ C\in [\k ]^\k :y\in C\}$, which implies that $\k \setminus \{ y \} \in \P (\r )\setminus\P (\s )$.
Thus $\r \not\sim _4 \s$.

Further we prove that $\sim _4 \;\not\subset \;\sim _3$ and, hence, $\sim _1 \;\varsubsetneq \;\sim _4$ and $\sim _3 \;\varsubsetneq \;\sim _6$.
Let $x\in \k$ and $\r :=\la \{ x\},\emptyset, \emptyset ,\dots  \ra $ and $\s :=\la \k \setminus \{ x \} ,\emptyset, \emptyset ,\dots  \ra$.
Then, clearly, $\r \not\cong \s$, that is $\r \not\sim _3 \s$. As above we have
$\P ( \r )=\{ C\in [\k ]^\k :x\in C\}$ and, by Lemma \ref{T6023}(a),
$\P ( \s )= \{ C_1 \cup C_2 : C_1 \in [\k \setminus \{ x \}]^\k \land C_2 \in [\{ x \}]^1 \}=\{ C\in [\k ]^\k :x\in C\}=\P ( \r )$.
Thus $\r \sim _4 \s$.

Finally we prove that $\sim _8 \;\neq \;\sim _{11}$, which implies $\sim _6 \;\neq \;\sim _{11}$.
Let $U\subset \k$, where $|U|=|\k \setminus U|=\k$ and let $\r :=\la \emptyset,\emptyset, \emptyset ,\dots  \ra $
and $\s :=\la U ,\emptyset, \emptyset ,\dots  \ra$. Then, by Lemma \ref{T6023} (b) and (c), $\P ( \r)=\la [\k ]^\k ,\subset \ra$,
$\P ( \s )=\la [\k ]^\k ,\subset \ra \times \la [\k ]^\k ,\subset \ra$, and
\begin{equation}\label{EQ6044}\textstyle
\sq \P ( \r) \cong (P(\k )/[\k ]^{<\k })^+ ,
\end{equation}
\begin{equation}\label{EQ6045}\textstyle
\sq \P ( \s) \cong (P(\k )/[\k ]^{<\k })^+ \times (P(\k )/[\k ]^{<\k })^+ .
\end{equation}
By Fact \ref{T6025}(d), the poset $(P(\k )/[\k ]^{<\k })^+$ is not isomorphic to its square. So, by (\ref{EQ6044}) and (\ref{EQ6045}) we have $\r \not\sim _8\s$.

(c) For $\r$ and $\s$ defined in the previous paragraph we have $\r \not\sim _8\s$ and we prove that $\r \sim _{10}\s$. First we consider the case when $\k >\o$.
By (\ref{EQ6044}) and Fact \ref{T6025}, $ \P ( \r) \equiv (\Col (\o , 2^\k ))^+$.
By (\ref{EQ6045}), forcing by the poset $\sq \P ( \s)$ collapses
$2^\k$ to $\o$ and, since the poset is of size $2^\k$, by Fact \ref{T6025}(c) we have $\ro \sq \P (\s) \cong \Col (\o , 2^\k )$.
Thus the posets $\P ( \r)$ and $\P ( \s)$ are forcing equivalent, that is $\r \sim _{10}\s$. If $\k =\o$ we use Fact \ref{T6025}(b).
\kdok
The following theorem shows that the equivalence of the similarities $\sim _8 $ (the isomorphism of $\sq \P(\X )$) and $\sim _{10}$
(the isomorphism of $\ro \sq \P(\X )$) is independent of ZFC even for the simplest unary language.
\begin{te} \label{T6028}
If $L$ is the language containing only one unary relational symbol, then on $\Mod _L (\o )$ we have
$\sim _8 \;=\; \sim _6$ and
$$
\sim _{10}\; =\left\{ \begin{array}{ll}
                                    \sim _{11} & \mbox{ if the poset }(P(\o )/\Fin )^+  \mbox{ is forcing equivalent to its square,}\\
                                    \sim _6    & \mbox{ otherwise.}
                   \end{array}
           \right.
$$
So, the equality $\sim _8 \;=\; \sim _{10}$ is independent of ZFC.
\end{te}
\dok
By Lemma \ref{T6023}, for $U\subset \o$, writing $\P (U)$ instead of $\la \P (\o ,U ), \subset \ra$, we have
\begin{equation}\label{EQ6046}
\P (U) \cong \left\{ \begin{array}{ll}
                                       \la [\o ]^\o , \subset \ra       &  \mbox{ if } |U|<\o \mbox{ or } |\o\setminus U|<\o ,\\
                                       \la [\o ]^\o , \subset \ra ^2  &  \mbox{ otherwise} ;
                      \end{array}
              \right.
\end{equation}
\begin{equation}\label{EQ6047}
\sq \P (U) \cong \left\{ \begin{array}{ll}
                                       (P(\o )/\Fin )^+ &  \mbox{ if } |U|<\o \mbox{ or } |\o\setminus U|<\o ,\\
                                       ((P(\o )/\Fin )^+ )^2   &  \mbox{ otherwise} .
                      \end{array}
              \right.
\end{equation}
If $U_1,U_2 \subset \o $ and $U_1 \not\sim _6 U_2$, that is $\P (U_1 ) \not\cong \P (U_2 )$, then, by
(\ref{EQ6046}) and (\ref{EQ6047}), for example, $\sq \P (U_1 )\cong (P(\o )/\Fin )^+$ and
$\sq  \P (U_2 ) \cong ((P(\o )/\Fin )^+ )^2$ and, by Fact \ref{T6025}(d),
$\sq \P (U_1 ) \not\cong \sq \P (U_2 )$, that is $U_1 \not\sim _8 U_2$.
Thus $\sim _8 \;\subset \; \sim _6$, which implies $\sim _8 \;= \; \sim _6$.

If $(P(\o )/\Fin )^+\equiv ((P(\o )/\Fin )^+)^2$, then by  (\ref{EQ6047}) for each $U\subset \o$ we have $\P (U)\equiv (P(\o )/\Fin )^+$ and,
hence $\sim _{10}=\sim _{11}$.
Otherwise, if $(P(\o )/\Fin )^+\not\equiv ((P(\o )/\Fin )^+)^2$, then for $U_1, U_2 \subset \o $ satisfying $U_1 \sim _{10} U_2$ by Fact
\ref{T4042}(b) we have $\sq \P (U_1 ) \equiv \sq \P (U_2 ) $ so,  by the assumption and (\ref{EQ6047}),
$\sq \P (U_1 ) \cong \sq \P (U_2 )$. Thus $\sim _{10}\;\subset \;\sim _8$ and, hence
$\sim _{10}\;= \;\sim _8\;= \;\sim _6$.

By Fact \ref{T6025}(b), CH implies that $(P(\o )/\Fin )^+\equiv ((P(\o )/\Fin )^+ )^2$.
But, by a result of Shelah and Spinas \cite{SheSpi}, in the Mathias model these two posets have different
distributivity numbers and, hence, they are not forcing equivalent.
\kdok

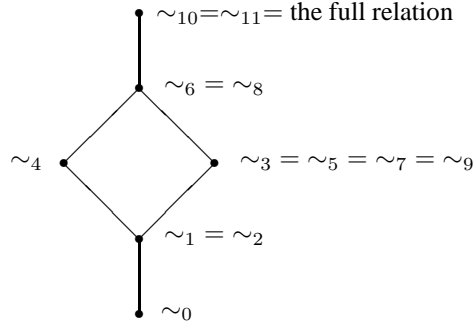
\begin{figure}[htbp]\label{FIG4}
\begin{center}
%TeXCAD Picture [konvergencije.pic]. Options:
%\grade{\on}
%\emlines{\off}
%\epic{\off}
%\beziermacro{\on}
%\reduce{\on}
%\snapping{\on}
%\pvinsert{\usepackage{amssymb,latexsym}}
%\quality{8.000}
%\graddiff{0.005}
%\snapasp{5}
%\zoom{2.0000}
\unitlength 1mm %0.7mm % = .854pt
\linethickness{0.5pt}
\ifx\plotpoint\undefined\newsavebox{\plotpoint}\fi % GNUPLOT compatibility

\vspace{-5mm}
\begin{picture}(80,50)(0,0)

%----------------------------- linije --------------------------

\put(30,5){\line(0,1){10}}%1
\put(30,15){\line(-1,1){10}}%2
\put(30,15){\line(1,1){10}}%3
\put(20,25){\line(1,1){10}}%4
\put(40,25){\line(-1,1){10}}%5
\put(30,35){\line(0,1){10}}%6

%----------------------------- tacke --------------------------

\put(30,5){\circle*{1}}%1
\put(30,15){\circle*{1}}%
\put(20,25){\circle*{1}}%
\put(40,25){\circle*{1}}%
\put(30,35){\circle*{1}}%
\put(30,45){\circle*{1}}%

%----------------------------- tekst -----------------------
\small
%\footnotesize
%\scriptsize
%\tiny
\put(40,15){\makebox(0,0)[cc]{$\sim _1 \;=\; \sim _2 $}}%
\put(15,25){\makebox(0,0)[cc]{$\sim _4 $}}%
\put(40,35){\makebox(0,0)[cc]{$\sim _6 \;=\; \sim _8$}}%
\put(59,25){\makebox(0,0)[cc]{$\sim _3 \;=\; \sim _5 \;=\; \sim _7 \;=\; \sim _9$}}%
\put(35,5){\makebox(0,0)[cc]{$\sim _0   $}}%1
\put(52,45){\makebox(0,0)[cc]{$\sim _{10} = \sim _{11} =$ the full relation}}%1
\end{picture}
\end{center}
\vspace{-7mm}
\caption{The similarities on $\Mod _{\la R \ra} (\o )$ if $(P(\o )/\Fin )^+\equiv ((P(\o )/\Fin )^+)^2$}
\end{figure}
%\vspace{-5mm}
\subsection{Infinite non-unary structures}
For infinite structures of non-unary languages the diagram from Figure \ref{F4008} does not collapse at all.
Namely the main result of this subsection is the following theorem.
\begin{te} \label{T6019}
If $L$ is a non-unary relational language and $\k$ an infinite cardinal, then in the diagram from Figure \ref{F4008} describing
the similarities $\sim _k$ on the set $\Mod _L (\k )$ all the implications $a$ - $o$ are proper and there are no new implications
(except the ones following from transitivity). Consequently, the same holds for the diagram from Figure \ref{F4002} related to the
class of all relational structures.
\end{te}
Theorem \ref{T6019} will be proved in two steps. First we will prove the statement for the class $\Mod _{L_b}(\o)$ of countable binary structures
(where $L_b =\la R \ra$ and $ar (R)=2$) and then, roughly speaking, make a correspondence between the classes $\Mod _{L_b}(\o)$ and $\Mod _L(\k )$
preserving all the similarities $\sim _k$ and their negations.

\subsubsection{Proof of Theorem \ref{T6019} for the class of countable binary structures}
First, giving examples (i.e.\ constructing pairs of structures), we show that for $L=L_b$ and $|X| =\o$, in the diagram from Figure \ref{F4008} all the implications $a$ - $o$ are proper. We will use the following auxiliary claim.
\begin{lem}      \label{T6006}
If $\P =\langle  P , \leq_P \rangle $ and $\Q =\langle  Q , \leq_Q \rangle $ are  partial orders and $f: P\rightarrow Q$ a surjection such that
for each $p_1 , p_2 \in P$ we have

(i) $p_1 \leq _P p_2 \Rightarrow f(p_1) \leq _Q^* f(p_2)$,

(ii) $p_1 \perp _P p_2 \Rightarrow f(p_1) \perp _Q f(p_2)$,

\noindent
then $\mathop{\rm sq}\nolimits  {\mathbb P} \cong \mathop{\rm sq}\nolimits  {\mathbb Q}$.
\end{lem}
\dok
First we prove that for each $p_1,p_2\in P$ we have
\begin{equation}\label{EQ6005}
p_1 \leq _P^* p_2 \Leftrightarrow f(p_1) \leq _Q^* f(p_2).
\end{equation}

($\Rightarrow$) Assuming $p_1 \leq _P^* p_2$ we have to prove that
\begin{equation}\label{EQ4093}
\forall q\leq _Q f(p_1) \;\;  q \not\perp _Q f(p_2) .
\end{equation}
Let $q\leq _Q f(p_1)$. Since $f$ is onto, there is $p_3\in P$ such that $f(p_3)=q$. Thus $f(p_3)\leq _Q f(p_1)$ and, by (ii), there is
$p_4 \leq _P p_3 , p_1$. So, since  $p_1 \leq _P^* p_2$ we have $p_4 \not\perp _P p_2$ namely there is $p_5\leq _P p_4, p_2$.
By (i) we have $f(p_5)\leq _Q^* f(p_2)$, which implies $f(p_5) \not\perp _Q f(p_2)$ and, hence, there is $q_0 \leq _Q f(p_5) ,f(p_2)$.
Since $p_5 \leq _P p_4 \leq _P p_3$, by (i) we have $f(p_5)\leq _Q^* f(p_3)=q$ and, hence, $q_0 \leq _Q^* q$, which implies $q_0 \not\perp _Q q$,
so there is $q' \leq _Q q_0, q$. Now $q' \leq _Q q, f(p_2)$ and (\ref{EQ4093}) is proved.

($\Leftarrow$) Assuming (\ref{EQ4093}) we prove that $p_1 \leq _P^* p_2$. So, taking $p\leq _P p_1$ we show that $p\not\perp _P p_2$.
By (i) we have $f(p)\leq _Q ^* f(p_1)$ which implies that there is $q \leq _Q f(p),f(p_1)$.
By (\ref{EQ4093}) we have $q \not\perp _Q f(p_2)$ and, hence, there is $q' \leq _Q q ,f(p_2 )$.
Now $q' \leq _Q f(p),f(p_2 )$ and, by (ii), $p\not\perp _P p_2$. Thus (\ref{EQ6005}) is proved.

Now we show that $\langle  P/\!\!=_P^* , \trianglelefteq _P \rangle \cong _F \langle  Q/\!\!=_Q^* , \trianglelefteq _Q \rangle $, where
$F([p])=[f(p)]$.
By (\ref{EQ6005}), for $p_1,p_2\in P$ we have
$[p_1]=[p_2]$
iff $p_1 =_P^* p_2$
iff $p_1 \leq_P^* p_2 \land p_2 \leq_P^* p_1$
iff $f(p_1) \leq_Q^* f(p_2) \land f(p_2) \leq_Q^* f(p_1)$
iff $f(p_1) =_Q^* f(p_2)$
iff $[f(p_1)]=[f(p_2)]$
iff $F([p_1])= F([p_2])$
and $F$ is a well defined injection.
Since $f$ is onto, for $q\in Q$ there is $p\in P$ such that $q=f(p) $. Thus $F([p])=[f(p)]= [q]$ and $F$ is onto.

By (\ref{EQ6005}) again,
$[p_1] \trianglelefteq _P [p_2]$
iff $p_1 \leq_P^* p_2$
iff $f(p_1) \leq_Q^* f(p_2) $
iff $[f(p_1)] \trianglelefteq _Q [f(p_2)]$
iff $F([p_1])\trianglelefteq _Q F([p_2])$.
Thus $F$ is an  isomorphism.
\kdok

\begin{ex}\rm\label{EX6008}
The implication $a$ can not be reversed.
Let $\X = \la \o , \leq \ra$ and $\Y =\la \o , \leq _f \ra$, where $f:\o\rightarrow \o$ is a bijection different from the identity and
$\leq _f =\{ \la f(m),f(n) \ra : m\leq n \}$. Then $\X \cong \Y$
and $\P (\X )=\P (\Y )=[\o ]^\o$, but $\X \neq \Y $.
\end{ex}
\begin{ex}\rm\label{EX6009}
The implications $b$ and $f$ can not be reversed. Let

$\X = \la \o , \{ \la n,n+1 \ra :n\in \o  \} \cup \{ \la 2n , 2n \ra : n\in \o  \}  \ra$ and

$\Y =\la \o , \{ \la n,n+1 \ra :n\in \o  \} \cup \{ \la 2n+1 , 2n+1 \ra : n\in \o  \}  \ra$.

\noindent
Then $\P (\X )=\P (\Y )=\{ [2n , \infty ) :n\in \o \}$
and $\X \rightleftarrows \Y$ but $\X \not\cong \Y $.
\end{ex}
\begin{ex}\rm\label{EX6001}
The implications $c$, $e$ and $g$ can not be reversed.
Let us define $\X = \la \o ,\o ^2 \setminus \{ \la 0,0 \ra\} \ra$ and $\Y =\la \o ,\o ^2 \setminus \{ \la 1,1 \ra\}  \ra$. Then $\X \cong \Y$ and
 $\P (\X )= \{ A\in [\o ]^\o : 0\in A \} \cong \P (\Y )=\{ A\in [\o ]^\o : 1\in A \}$, but $\P (\X )\neq \P (\Y )$.
\end{ex}
\begin{ex}\rm\label{EX6000}
The implications $d$, $h$, $k$ and $n$ can not be reversed.
Let $\X = \la \o , \leq \ra$ and $\Y =\la \o , \o \times \o \ra$. Then $\P (\X )=\P (\Y )=[\o ]^\o$ and, hence, $\P (\X )\cong \P (\Y )$,
$\sq \P (\X )\cong \sq \P (\Y )$ and $\P (\X )\equiv \P (\Y )$, but $\X \not\rightleftarrows \Y$.
\end{ex}
\begin{ex}\rm\label{EX6010}
The implications $i$ and $j$ can not be reversed.
Let $\X = \la (0,1)_\Q , \leq \ra$ and $\Y =\la (0,1]_\Q , \leq \ra$ be suborders of the rational line, $\Q$.
Then, clearly, $\X \rightleftarrows\Y$.

Since the elements of $\P (\X )$ are dense linear orders without end points,
 each chain $\L$ in the poset $\la \P (\X ), \subset\ra$ has a supremum: $\bigcup \L$.  On the other hand,
$\L =\{ (0, \frac{1}{2}-\frac{1}{n} ]_\Q : n\geq 3\}$ is a chain in the poset  $\la \P (\Y ), \subset\ra$,
$\bigcup \L =(0,\frac{1}{2})\not\in \P (\Y )$ and the sets $(0,\frac{1}{2})_\Q \cup \{ q \}$, $q\in [\frac{1}{2}, 1]_\Q$,
are upper bounds for $\L$, but $\L$ does not have a least upper bound. Thus the poset  $\la \P (\Y ), \subset\ra$ is not chain complete and, hence,
$\P (\X )\not\cong \P (\Y )$.

Using Lemma \ref{T6006} we show that $\sq \P (\X )\cong \sq \P (\Y )$.
We remind the reader that a linear order $L$ is called scattered iff $\Q \not\hookrightarrow L$. Let $\Scatt$ denote the set of scattered
suborders of $\Q$. It is easy to see that for $A,B\in \P (\X )$ we have $A\leq ^*B \Leftrightarrow A\setminus B\in \Scatt$ and
$A\perp B \Leftrightarrow A\cap B\in \Scatt$ (where $\leq ^*$ is the corresponding separative modification) and that the same holds
for $A,B\in \P (\Y )$. Clearly, if $A\in \P (\Y)$, then $A\setminus \{ \max A \} \subset (0,1)_\Q$ and it is a copy of $\X$, so, the function
$f:\P (\Y )\rightarrow \P (\X )$, given by $f(A)=A\setminus \{ \max A \}$, is well defined and we show that it satisfies the assumptions of Lemma \ref{T6006}.
First, if $C\in \P (\X )$, then $C\subset (0,1)_\Q$ and, clearly, $C\cup \{ 1\} \in \P (\Y )$ and $f(C\cup \{ 1\})=C$. Thus $f$ is a surjection.
Let $A,B\in \P (\Y )$. If $A\subset B$, then $f(A)\setminus f(B)= (A\setminus \{ \max A \}) \setminus (B\setminus \{ \max B \})\subset \{ \max B \}\in \Scatt$ and, hence, $f(A)\leq ^* f(B)$ so (i) is true. If $A\perp B$, that is $A\cap B\in \Scatt$, then, clearly, $f(A)\cap f(B)\in \Scatt$, thus $f(A)\perp f(B)$ and
(ii) is true as well. By Lemma \ref{T6006} we have $\sq \P (\X )\cong \sq \P (\Y )$.
\end{ex}
\begin{ex}\rm\label{EX6004}
The implication $m$ can not be reversed.
By Example 4.4 of \cite{Ktow}, if $\X$ is the directed graph $\la {}^{<\o }2 , \r\ra$, where
$\r =\{ \la \f , \f ^\smallfrown i \ra : \f \in {}^{<\o }2 \land i\in 2 \}$, then
$\la \P (\X ),\subset \ra =\sq \la \P (\X ),\subset \ra\cong \la {}^{<\o }2 , \supset \ra$.
Let $\Y$ be the directed graph $\la {}^{<\o }3 , \s \ra$, where
$\s =\{ \la \f , \f ^\smallfrown i \ra : \f \in {}^{<\o }3 \land i\in 3 \}$, then in a similar way we show that
$\la \P (\Y ),\subset \ra =\sq \la \P (\Y ),\subset \ra\cong \la {}^{<\o }3 , \supset \ra$.
Clearly $\sq \la \P (\X ),\subset \ra \not\cong \sq \la \P (\Y ),\subset \ra$, but
$\ro \sq \la \P (\X ),\subset \ra \cong \ro \sq \la \P (\Y ),\subset \ra \cong \Borel /\M$.
\end{ex}
\begin{ex}\rm\label{EX6006}
The implication $l$ can not be reversed.
Let $\X$ be the directed graph from Example \ref{EX6004} and let $\Y$ be the directed graph $\la Y, \s\ra$,
where $Y\subset {}^{<\o }2$ and $\s \subset Y\times Y$ are defined by
$$
Y=\{ \emptyset ,0,1 \} \cup \{ jj^\smallfrown \f :j\in 2 \land \f \in {}^{<\o }2 \} ,
$$
$$
\s= \{ \la \emptyset , 0 \ra , \la \emptyset , 1 \ra , \la 0,00 \ra , \la 1, 11 \ra \} \cup
    \{ \la jj^\smallfrown \f , jj^\smallfrown \f ^\smallfrown k \ra : j,k\in 2 \land \f \in {}^{<\o }2  \} .
$$
It is  easy to see that $\X \rightleftarrows \Y $ and
$$
\P (\Y )= \{ Y \} \cup \{ A_{jj^\smallfrown \f }^{kl}: j,k,l \in 2 \land \f \in {}^{<\o }2\} ,
$$
where
$
A_{jj^\smallfrown \f }^{kl}
=\{ jj^\smallfrown \f ,jj^\smallfrown \f ^\smallfrown 0 , jj^\smallfrown \f ^\smallfrown 1\} \cup
\{ jj^\smallfrown \f ^\smallfrown 0 ^\smallfrown k ^\smallfrown \p : \p \in {}^{<\o }2 \} \cup
\{ jj^\smallfrown \f ^\smallfrown 1 ^\smallfrown l ^\smallfrown \p : \p \in {}^{<\o }2 \} .
$
By Example \ref{EX6004}, the poset $ \sq \la \P (\X ),\subset \ra$
is isomorphic to the reversed binary tree. Thus, in order to prove that $\sq \P (\Y )\not\cong \sq \P (\X )$
we will show that $[A_{00}^{00}]$ and $[A^{01}_{00}]$ are incomparable but compatible elements of
$\sq \P (\Y )=\la \P (\Y )/\!=^* , \trianglelefteq \ra$. So we have
$$
A_{00}^{00}
=\{ 00 ,000 , 001\} \cup
\{ 0000 ^\smallfrown \p : \p \in {}^{<\o }2 \} \cup \{ 0010 ^\smallfrown \p : \p \in {}^{<\o }2 \} ,
$$
$$
A_{00}^{01}
=\{ 00 ,000 , 001\} \cup
\{ 0000 ^\smallfrown \p : \p \in {}^{<\o }2 \} \cup \{ 0011 ^\smallfrown \p : \p \in {}^{<\o }2 \} .
$$
Clearly $\{ 0000 ^\smallfrown \p : \p \in {}^{<\o }2 \}$ is a copy of $\X $ and, hence, contains a copy of $\Y$, say $B$.
Since $B\subset A_{00}^{00}, A_{00}^{01}$ we have $B \leq ^* A_{00}^{00}, A_{00}^{01}$ and $[B] \trianglelefteq [A_{00}^{00}], [A_{00}^{01}]$ thus
 $[A_{00}^{00}]$ and $[A^{01}_{00}]$ are compatible elements of $\sq \P (\Y )$.

In order to prove that  $[A_{00}^{00}]\not\trianglelefteq [A_{00}^{01}]$ we need $C\in \P (\Y )$
such that $C\subset A_{00}^{00}$ and $D\not\subset C\cap A_{00}^{01} $, for all $D\in \P (\Y )$.
Now $\{ 0010 ^\smallfrown \p : \p \in {}^{<\o }2 \}\subset A_{00}^{00}$ is a copy of $\X $ and, hence, contains a copy of $\Y$, say $C$.
Since $\{ 0010 ^\smallfrown \p : \p \in {}^{<\o }2 \}\cap A_{00}^{01}=\emptyset$, we have $C \cap A_{00}^{01}=\emptyset$ and we are done.
Thus $[A_{00}^{00}]\not\trianglelefteq [A_{00}^{01}]$ and, similarly, $[A_{00}^{01}]\not\trianglelefteq [A_{00}^{00}]$.
\end{ex}
Thus in Figure \ref{F4008} for $\Mod _{L_b}(\o )$ all the implications $a$ - $o$ are proper and we show that there are no new implications
except the ones following from transitivity. So it remains to be shown that the eight pairs which are incomparable in the Hasse diagram  in Figure \ref{F4008}
are really incomparable. We will use the following elementary fact: if $\P =\la P ,\leq \ra$ is a partial order and $p,q,r\in P$, then
\begin{equation}\label{EQ6050}
r=p\wedge q \;\mbox{ and }\;  r< p \;\mbox{ and }\;  r< q \;\Rightarrow \; p\parallel q.
\end{equation}
In fact our poset of similarities is a
suborder of the lattice  $\la EQ (\Int _{L_b} (\o )), \subset \ra $ of equivalence relations on the set $\Int _{L_b} (\o)$, where
for  $\sim , \sim ' \in EQ (\Int _{L_b} (\o ))$ we have $\sim \wedge \sim ' =\sim \cap \sim '$ and $\sim \vee \sim ' =trcl (\sim \cup \sim ')$
and $\sim \subset \sim ' $ iff the $\sim$-partition is a refinement of the $\sim '$-partition of $\Int _{L_b} (\o )$. Now, since {\it by our definition}
we have $\sim _1 = \sim _2 \cap \sim _3$, by (\ref{EQ6050}) we obtain $\sim _2 \;\parallel\; \sim _3$ and similarly for the other seven pairs.
\subsubsection{Proof of Theorem \ref{T6019}}
The following concepts and facts will be used in our proof. Let $L_{b}=\la R\ra$, where $\ar (R)=2$. If $\X=\la X,\r \ra$ is an $L_b$-structure, then
the transitive closure $\r _{rst}$ of the relation $\r _{rs} =\Delta _X \cup \r \cup \r ^{-1}$ (given by $x \,\r_{rst} \,y$ iff there are $n\in \N$ and
$z_0 =x , z_1, \dots ,z_n =y$ such that $z_i \;\r _{rs} \;z_{i+1}$, for each $i<n$)
is the minimal equivalence relation on $X$ containing $\r$. The corresponding equivalence classes are called the {\it components} of $\X$ and
the structure $\X$ is called {\it connected} iff $|X/\r _{rst} |=1$. The {\it complement} of the structure $\X $,
$\la X, (X\times X)\setminus \r \ra $ will be denoted by $\X ^c$; its {\it reflexification}, $\la X, \r \cup \Delta _X \ra$,  by $\X _{re} $; and
its {\it irreflexification}, $\la X, \r \setminus \Delta _X \ra$, by $\X _{ir} $.

If $\X _i=\la X_i, \r _i \ra$, $i\in I$, are connected $L_b$-structures  and $X_i \cap X_j =\emptyset$, for
different $i,j\in I$, then the structure $\bigcup _{i\in I} \X _i =\la \bigcup _{i\in I} X_i , \bigcup _{i\in I} \r _i\ra$ is
the {\it disjoint union} of the structures $\X _i$, $i\in I$,  and the structures $\X _i$, $i\in I$, are its components.
\begin{fac}\label{T4010}
(\cite{Ktow})
If $\X$ is an $L_b$-structure, then at least one of the structures $\X$ and $\X ^c $ is connected.
\end{fac}
\begin{fac}\label{T4015}(\cite{Ktow})
Let $\X _i = \la X_i , \r _i \ra , i\in I$, and $\Y _j = \la Y_j , \s _j \ra , j\in J$, be families of disjoint connected
binary structures.
Then $F:\bigcup _{i\in I}\X _i  \hookrightarrow \bigcup _{j\in J}\Y _j$ iff there are $f:I\rightarrow J$ and $g_i : \X _i \hookrightarrow \Y _{f(i)}$, $i\in I$,
such that $F=\bigcup _{i\in I} g_i$ and
$\la g_i(x), g_{i'}(x')\ra \not\in\s _{rs}$, whenever $i\neq i'$, $x \in X_i$ and $x' \in X_{i'}$.
\end{fac}
\begin{fac}\label{T6102}
Let $\X $ be a binary structure. Then

(a) $\Emb (\X )= \Emb (\X  ^c )$ and $\P (\X )= \P (\X  ^c )$;

(b) If $\X $ is irreflexive, then $\Emb (\X )= \Emb (\X _{re})$ and $\P (\X )= \P (\X _{re})$;

(c) If $\X $ is reflexive, then $\Emb (\X )= \Emb (\X _{ir})$ and $\P (\X )= \P (\X _{ir})$.
\end{fac}
\begin{te}\label{T6017}(Vop\v enka, Pultr, Hedrl\'in \cite{Vop})
On any set $X$ there is an irreflexive binary relation $\r$ such that $\id _X$ is the only endomorphism of the structure $\la X,\r \ra$.
\end{te}
For a cardinal $\l$ let
$\Int ^*_{L_{b}}(\l ) = \{ \r \subset \l ^2 :  \la \l , \r \ra \mbox{ is connected } \land \; \r \cap \Delta _\l \neq \emptyset \} $.
Then $\Int ^*_{L_{b}}(\l ) \subset \Int _{L_b}(\l ) $ and
$\Mod ^*_{L_{b}}(\l ) := \{ \la \l , \r \ra :\r \in \Int ^*_{L_{b}}(\l )  \}\subset \Mod _{L_b}(\l ) $.

\begin{te}\label{T6018}
Let $\k \geq \l \geq \o$ be cardinals and $L=\la R_i : i\in I\ra$ a non-unary relational language. Then there is a mapping
$\t : \Int ^*_{L_{b}}(\l ) \rightarrow \Int _L (\k )$ such that

(a) $\P (\k , \t _\r )\cong \P (\l ,\r)  $, for each $\r \in \Int ^*_{L_{b}}(\l )$;

(b) For each $\r \in \Int _{L_b}(\l )$ there are $\r '\in \Int ^*_{L_{b}}(\l )$ and $\t \in \Int _L(\k )$ such that
$\P (\l , \r ')=\P (\l , \r )\cong \P (\k , \t )$;

(c) $\tau$ preserves all the relations $\sim _k$ from Figure \ref{F4008}, that is for each $k \leq 11$
\begin{equation}\label{EQ6021}
\forall \r ,\s \in \Int ^*_{L_b}(\l ) \;\; (\r \sim _k \s \Leftrightarrow \tau _\r \sim _k \tau _\s ).
\end{equation}
\end{te}
\dok
First suppose that $\l <\k$. Then $|\k \setminus \l|=\k$ and, by Theorem \ref{T6017} we can fix an irreflexive
binary relation $\theta \subset (\k \setminus \l)^2 $ such that $\Emb (\k \setminus \l ,\theta )=\{ \id _{\k \setminus \l }\}$.
By Theorem \ref{T4010} and Facts \ref{T6102}(a) and \ref{T6102}(c) we can assume that the relation $\theta $ is connected and irreflexive.
The language $L$ is not unary and we fix an $i_0\in I$ such that $n_{i_0}\geq 2$.
Now, for $\r \in \Int ^*_{L_{b}}(\l )$ let the interpretation $\t _\r =\la \t ^\r_i : i\in I\ra \in \Int _L (\k )$ be defined by
\begin{equation}\label{EQ6022}
\t ^\r _i = \left\{ \begin{array}{cl}
                                       (\r \cup \theta)\times \k ^{n_{i_0}-2} &  \mbox{ if } i=i_0 \mbox{ and } n_{i_0}>2;\\
                                       (\r \cup \theta)                       &  \mbox{ if } i=i_0 \mbox{ and } n_{i_0}=2;\\
                                       \emptyset                              &  \mbox{ if } i\neq i_0 .
                      \end{array}
              \right.
\end{equation}
For convenience, for  $\r ,\s \in \Int ^*_{L_{b}}(\l )$, instead of
$\Emb (\la\k , \t _\r \ra , \la \k , \t _\s \ra)$ (respectively, $\Emb (\la\l , \r \ra , \la \l , \s \ra)$) we will write
$\Emb (\t _\r , \t _\s )$ (resp.\ $\Emb ( \r  ,  \s )$).
\begin{cla}\label{T6033}
For each $\r ,\s \in \Int ^*_{L_{b}}(\l )$  we have

(i) $\Emb ( \t _\r ,  \t _\s )=\{ f\cup \id _{\k \setminus \l } : f\in \Emb ( \r , \s )\}$;

(ii) $\Iso (\t _\r ,  \t _\s)=\{ f\cup \id _{\k \setminus \l } : f\in \Iso (\r , \s )\}$;

(iii) $\P (\t _\r ,  \t _\s)=\{ C\cup (\k \setminus \l ) : C\in \P (\r , \s )\}$;

(iv) $\Emb ( \t _\r  )=\{ f\cup \id _{\k \setminus \l } : f\in \Emb ( \r  )\}$;

(v) $\Aut ( \t _\r  )=\{ f\cup \id _{\k \setminus \l } : f\in \Aut ( \r  )\}$;

(vi) $\P ( \t _\r  )=\{ C\cup (\k \setminus \l ) : C\in \P ( \r  )\}$.
\end{cla}
\dok
For convenience let $\pi_\r := \r \cup \theta$, for $\r  \in \Int ^*_{L_{b}}(\l )$. First we prove that
\begin{equation}\label{EQ6023}
\Emb (\la \k, \pi _\r \ra , \la \k, \pi _\s  \ra ) =\{ f\cup \id _{\k \setminus \l } : f\in \Emb ( \r  , \s )\} .
\end{equation}
By the construction, $\la \k, \pi _\r \ra =\la \l ,\r \ra \cup \la \k \setminus \l, \theta\ra$ and
$\la \k, \pi _\s \ra =\la \l ,\s \ra \cup \la \k \setminus \l, \theta\ra$ are partitions of the binary structures
$\la \k, \pi _\r \ra$ and $\la \k, \pi _\s \ra$ into their connectivity components.
Since $\r \cap \Delta _\l \neq \emptyset$  and $\theta$ is an irreflexive relation,
we have $\la \l ,\r \ra\not\hookrightarrow \la \k \setminus \l, \theta\ra$ and the inequality $\k >\l$ implies that
$\la \k \setminus \l, \theta\ra \not\hookrightarrow \la \l ,\s \ra$. So, by Theorem \ref{T4015},
$F\in \Emb (\la \k, \pi _\r \ra , \la \k, \pi _\s  \ra )$ iff $F\upharpoonright \l \in \Emb (\la\l , \r \ra , \la \l , \s \ra)$ and
$F\upharpoonright (\k \setminus\l) \in \Emb (\la \k \setminus \l, \theta\ra )=\{ \id _{\k \setminus \l }\}$ and (\ref{EQ6023}) is proved.

Now we prove
\begin{equation}\label{EQ6024}
\Emb (\la \k, \t ^\r_{i_0} \ra , \la \k, \t ^\s _{i_0}  \ra ) =\{ f\cup \id _{\k \setminus \l } : f\in \Emb (\r  , \s )\} .
\end{equation}
If $F:\k \rightarrow \k$ is an injection, then
$F\in \Emb (\la \k, \t ^\r_{i_0} \ra , \la \k, \t ^\s _{i_0}  \ra )$ iff
for each $x_1, x_2 , \dots , x_{n_{i_0}} \in \k $
$$
\la x_1, x_2 , \dots , x_{n_{i_0}}\ra \in \pi _\r \times \k ^{n_{i_0}-2} \Leftrightarrow
\la F(x_1), F(x_2) , \dots , F(x_{n_{i_0}})\ra \in \pi _\s \times \k ^{n_{i_0}-2}
$$
iff for each $x_1, x_2 \in \k $ we have: $\la x_1, x_2 \ra \in \pi _\r  \Leftrightarrow \la F(x_1), F(x_2) \ra \in \pi _\s $, iff
$F\in \Emb (\la \k, \pi _\r \ra , \la \k, \pi _\s  \ra )$. Now (\ref{EQ6024}) follows from (\ref{EQ6023}).

(i) Clearly, $F\in \Emb (\t _\r , \t _\s )$ iff $F\in \Emb (\la\k , \t ^\r _i \ra , \la \k , \t ^\s_i \ra)$, for all $i\in I$.
By (\ref{EQ6022}) this holds iff $F\in \Emb (\la\k , \t ^\r_{i_0} \ra , \la \k ,\t ^\s _{i_0}  \ra)$ and we apply (\ref{EQ6024}).

(ii) If $f\in \Emb ( \r  ,  \s )$ then $f\cup \id _{\k \setminus \l }$ is a surjection iff $f$ is a surjection iff
$f\in \Iso ( \r  ,  \s )$. Now we apply (i).

(iii) $A\in \P (\t _\r ,  \t _\s )$ iff there is $F\in \Emb ( \t _\r , \t _\s )$ such that
$A=F[\k ]$ so, by (i), iff $A=f[\l ]\cup (\k \setminus \l)$, for some $f\in \Emb ( \r ,  \s )$, iff
$A=C\cup (\k \setminus \l)$, for some $C\in \P ( \r ,  \s )$.

Statements (iv), (v) and (vi) follow from (i), (ii) and (iii) respectively.
\kdok

Now we prove the theorem.

(a) By Claim \ref{T6033}(vi) we have $\P (\t _\r )=\{ C\cup (\k \setminus \l ) : C\in \P ( \r  )\}$ and
it is easy to check that the mapping $F: \P (\r) \rightarrow \P (\t _\r )$, defined by $F(C)=C\cup (\k \setminus \l )$,
is an isomorphism of the posets  $\la \P (\r) , \subset \ra$ and $\la \P ( \t _\r ) ,\subset \ra$.

(b) Let $\r\in \Int _{L_b}(\l )\setminus \Int ^*_{L_{b}}(\l )$. If $\r$ is connected, then  it is irreflexive, thus
$\r _{re} \in \Int ^*_{L_{b}}(\l )$ and, by Fact \ref{T6102}(b), $\P (\l , \r _{re})=\P (\l , \r )$.
Otherwise, by Theorem \ref{T4010} the relation $\r ^c$ is connected and, by Fact \ref{T6102}(a), $\P (\l , \r ^c)=\P (\l , \r )$.
Now, if $\r ^c \cap \Delta _\l \neq \emptyset$, we have $\r ^c \in \Int ^*_{L_{b}}(\l )$;
otherwise  $(\r ^c)_{re} \in \Int ^*_{L_{b}}(\l )$ and, by Fact \ref{T6102}(b),
$\P (\l , (\r ^c)_{re})=\P (\l , \r ^c )=\P (\l , \r )$.

If $\r \in \Int _{L_b}(\l )$ and $\r ' \in \Int ^*_{L_b}(\l )$, where $\P (\l ,\r )= \P (\l ,\r ')$, then
by (a) we have $\P (\l ,\r ')\cong \P (\k , \t _{\r '})$, where $\t _{\r '}\in \Int _L (\k )$. Thus $\P (\k , \t _{\r '})\cong\P (\l ,\r ) $.

(c) It is sufficient to prove that the mapping $\t : \Int ^*_{L_b}(\l ) \rightarrow \Int _L (\k )$
preserves the relations $\sim _k$, for $k\in \{0,3,4,6,8,9,10\}$.
Let $\r ,\s \in \Int ^*_{L_b}(\l )$.

$\sim _0$: $\r =\s \Leftrightarrow \t _\r =\t _\s$.
By (\ref{EQ6022}) we have: $\t _\r =\t _\s$ iff $\t ^\r _{i_0} =\t ^\r _{i_0}$ iff
$\r \cup \theta =\s \cup \theta$ iff $\r = \s$.

$\sim _3$: $\r \cong \s \Leftrightarrow \t _\r \cong \t _\s$.
If $\r \cong \s$ and $f\in \Iso ( \r , \s )$, then, by Claim \ref{T6033}(ii),
$f\cup \id _{\k \setminus \l }\in \Iso ( \t _\r ,  \t _\s )$ and, hence, $\t _\r \cong \t _\s$. Conversely, if $\t _\r \cong \t _\s$ and
$F\in \Iso ( \t _\r ,\t _\s )$, then, by Claim \ref{T6033}(ii), $F\upharpoonright \l \in \Iso ( \r , \s )$ and, hence, $\r \cong \s$.

$\sim _9$: $\r \rightleftarrows \s \Leftrightarrow \t _\r \rightleftarrows \t _\s$.
If $\r \hookrightarrow \s$ and
$f\in \Emb ( \r , \s )$, then, by Claim \ref{T6033}(i),
$f\cup \id _{\k \setminus \l }\in \Emb ( \t _\r , \t _\s )$ and, hence, $\t _\r \hookrightarrow \t _\s$. Thus
$\r \rightleftarrows \s $ implies $\t _\r \rightleftarrows \t _\s$.
Conversely, if $\t _\r \hookrightarrow \t _\s$ and
$F\in \Emb ( \t _\r ,  \t _\s )$, then, by Claim \ref{T6033}(i), $F\upharpoonright \l \in \Emb ( \r , \s )$ and, hence,
$\r \hookrightarrow \s$. So $\t _\r \rightleftarrows \t _\s$ implies $\r \rightleftarrows \s $.

$\sim _4$: $\P (\r ) =\P (\s )\Leftrightarrow \P (\t _\r )=\P (\t _\s )$. This follows from Claim \ref{T6033}(vi).

$\sim _6$: $\P (\r ) \cong \P (\s )\Leftrightarrow \P (\t _\r )\cong \P (\t _\s )$. This is true since by (a) we have
\begin{equation}\label{EQ6026}
\P (\r ) \cong \P (\t _\r ) \;\; \mbox{ and }\;\; \P (\s ) \cong \P (\t _\s ).
\end{equation}

$\sim _8$: $\sq \P (\r ) \cong \sq \P (\s )\Leftrightarrow \sq \P (\t _\r )\cong \sq \P ( \t _\s )$. This is true since by (\ref{EQ6026}) and Fact \ref{T4042}(a)
we have $\sq \P (\r ) \cong \sq \P (\t _\r )$ and $\sq \P (\s ) \cong \sq \P (\t _\s )$.

$\sim _{10}$: $\ro \sq \P (\r ) \cong \ro \sq \P (\s )\Leftrightarrow \ro \sq \P (\t _\r )\cong \ro \sq \P (\t _\s )$.
By (\ref{EQ6026}) and Fact \ref{T4042}(a)
we have $\ro \sq \P (\r ) \cong \ro \sq \P (\t _\r )$ and $\ro \sq \P (\s ) \cong \ro \sq \P (\t _\s )$.

So, the theorem is proved for $\l<\k$.  If $\l =\k$, then we define $\t ^\r _{i_0}:= \r \times \k ^{n_{i_0}-2}$
and continue in the same way.
\kdok
Finally we prove Theorem \ref{T6019}. In Subsection 3.2.1 it is shown that all the implications $a$ - $o$ in Figure \ref{F4008} for
the class $\Mod _{L_b}(\o)$ are proper. For example, concerning the implication $a$, in Example \ref{EX6008} we have constructed $\r ,\s \in \Int ^* _{L_b}(\o )$
such that $\r\sim _1\s$ but $\r\not\sim _0\s$. By Theorem \ref{T6018}(c) we have $\t _\r \sim _1\t _\s$ and $\t _\r\not\sim _0\t _\s$, which implies that
in Figure \ref{F4008} for the class $\Mod _L(\k )$ the implication $a$ is proper as well. The reader will notice that the structures constructed in
Examples \ref{EX6008} - \ref{EX6010}  belong to $\Int ^* _{L_b}(\o )$ and that the structures constructed in
Examples \ref{EX6004} and \ref{EX6006} are irreflexive. But their refexifications are in $\Int ^* _{L_b}(\o )$. Thus all the implications $a$ - $o$ in Figure \ref{F4008} for the class $\Mod _L(\k)$ are proper and using the same argument as in Subsection 3.2.1 we conclude that there are no additional implications
in the diagram describing the hierarchy of the considered similarities on the class $\Mod _L(\k)$.


\footnotesize
\begin{thebibliography}{99}
\bibitem{BalcVop}
      B.\ Balcar, P.\ Vop\v enka,
      On systems of almost disjoint sets,
      Bull.\ Acad.\ Polon.\ Sci.\ S\'er.\ Sci.\ Math.\ Astronom.\ Phys.\ 20 (1972) 421--424.
\bibitem{BalcSim}
      B.\ Balcar, P.\ Simon,
      Disjoint refinement,
      in: J.\ D.\ Monk and R.\ Bonnet (Eds.), Handbook of Boolean algebras, Vol.\ 2, 333--388,
      Elsevier Science Publishers B.V., Amsterdam, 1989.
\bibitem{Jech}
      T.\ Jech,
      Set Theory, 2nd corr.\ Edition,
      Springer, Berlin, 1997.
\bibitem{Kopp1}
      S.\ Koppelberg,
      General Theory of Boolean Algebras,
      in: J.\ D.\ Monk and R.\ Bonnet (Eds.), Handbook of Boolean Algebras, Part I,
      Elsevier Science Publishers B.V., Amsterdam, 1989.
\bibitem{Kun}
      K.\ Kunen,
      Set Theory, An Introduction to Independence Proofs,
      North-Holland, Amsterdam, 1980.
\bibitem{Ktow}
      M.\ S.\ Kurili\'c,
      From $A_1$ to $D_5$: Towards a forcing-related classification of relational structures,
      J.\ Symbolic Logic 79,1 (2014) 279--295.
\bibitem{Kemb}
      M.\ S. Kurili\'c,
      Maximally embeddable components,
      Arch.\ Math.\ Logic 52,7 (2013) 793--808.
\bibitem{Kurord}
      M.\ S. Kurili\'c,
      Forcing with copies of countable ordinals,
      Proc.\ Amer.\ Math.\ Soc.\ (in print).
\bibitem{Kur1}
      M.\ S. Kurili\'c,
      Posets of copies of countable scattered linear orders,
      Ann.\ Pure  Appl.\ Logic 165 (2014) 895--912.
\bibitem{Kstr}
      M.\, S. Kurili\'c,
      Isomorphic and strongly connected components,
      Arch.\ Math.\ Logic (in print) DOI10.1007/s00153-014-0399-2.
\bibitem{KurTod}
      M.\ S. Kurili\'c, S. Todor\v cevi\'c,
      Forcing by non-scattered sets,
      Ann.\ Pure Appl.\ Logic 163 (2012) 1299--1308.
\bibitem{KurTod1}
      M.\ S. Kurili\'c, S. Todor\v cevi\'c,
      Copies of the random graph,
      (to appear).
\bibitem{KurTod2}
      M.\ S. Kurili\'c, S. Todor\v cevi\'c,
      Copies of the random graph: the 2-localization,
      (to appear).
\bibitem{SheSpi}
      S.\ Shelah, O.\ Spinas,
      The distributivity numbers of $P(\omega)/$fin and its square,
      Trans.\ Amer.\ Math.\ Soc.\  352,5  (2000) 2023--2047.
\bibitem{Vop}
      P.\ Vop\v enka, A.\ Pultr, Z.\ Hedrl\'in,
      A rigid relation exists on any set,
      Comment.\ Math.\ Univ.\ Carolinae 6 (1965) 149--155.
\end{thebibliography}
\end{document}